\documentclass[final]{siamltex}

\usepackage[tbtags]{amsmath}
\usepackage{amssymb}
\usepackage[misc]{ifsym}
\usepackage{colortbl}
\numberwithin{equation}{section}   

\title{Connection between MP and DPP for Stochastic Recursive Optimal Control Problems: Viscosity Solution Framework in General Case}

\author{Tianyang Nie\thanks{School of Mathematics, Shandong University, Jinan 250100, P.R.China. E-mail: nietianyang@sdu.edu.cn. This author thanks the financial support from the National Natural Sciences Foundations of China (11601285), the Natural Science Foundation of Shandong Province (ZR2016AQ13, ZR2015JL003), and the Fundamental Research Fund of Shandong University (2015HW023).}\and Jingtao Shi\thanks{Corresponding author. School of Mathematics, Shandong University, Jinan 250100, P.R.China. E-mail: shijingtao@sdu.edu.cn. This author thanks the financial support from the National Natural Science Foundations of China (11301011, 11571205), and the Natural Science Fund for Distinguished Young Scholars of Shandong Province of China (JQ201401).}\and Zhen Wu\thanks{School of Mathematics, Shandong University, Jinan 250100, P.R.China. E-mail: wuzhen@sdu.edu.cn. This author thanks the financial support from the National Natural Science Foundation of China (61573217), 111 Project (B12023) and the Chang Jiang Scholar Program of Chinese Education Ministry.}}

\begin{document}
\maketitle

\begin{abstract}
This paper deals with a stochastic recursive optimal control problem, where the diffusion coefficient depends on the control variable and the control domain is not necessarily convex. We focus on the connection between the general maximum principle and the dynamic programming principle for such control problem without the assumption that the value is smooth enough, the set inclusions among the sub- and super-jets of the value function and the first-order and second-order adjoint processes as well as the generalized Hamiltonian function are established. Moreover, by comparing these results with the classical ones in Yong and Zhou [{\em Stochastic Controls: Hamiltonian Systems and HJB Equations, Springer-Verlag, New York, 1999}], it is natural to obtain the first- and second-order adjoint equations of Hu [{\em Direct method on stochastic maximum principle for optimization with recursive utilities, arXiv:1507.03567v1 [math.OC], 13 Jul. 2015}].
\end{abstract}



\pagestyle{myheadings}
\thispagestyle{plain}
\markboth{Connection between MP and DPP for SROCP}{Tianyang NIE, Jingtao SHI, Zhen WU}

\section{Introduction}

There are usually two ways to study optimal control problems: Pontryagin's {\it maximum principle} (MP) and Bellman's {\it dynamic programming principle} (DPP), involving the adjoint variable $\psi$, the Hamiltonian function $H$, and the value function $V$, respectively. The classical result on the connection between the MP and the DPP for the deterministic optimal control problem can be seen in Fleming and Rishel \cite{FR75}, which is known as $\psi(t)=-V_x(t,\bar{x}(t))$ and $V_t(t,\bar{x}(t))=H(t,\bar{x}(t),\bar{u}(t),\psi(t))$, where $\bar{u}$ is the optimal control and $\bar{x}$ is the optimal state. Since the value function $V$ is not always smooth, some non-smooth versions of the classical result were studied by using non-smooth analysis and generalized derivatives. An attempt to relate the MP and the DPP without assuming the smoothness of the value function was first made by Barron and Jessen \cite{BJ86}, where the viscosity solution was used to derive the MP from the DPP. Within the framework of viscosity solution, Zhou \cite{Zhou90} showed that
\begin{equation}\label{Zhou90}
\left\{
\begin{aligned}
&D_x^{1,-}V(t,\bar{x}(t))\subset\{-\psi(t)\}\subset D_x^{1,+}V(t,\bar{x}(t)),\\
&D_t^-V(t,\bar{x}(t))\subset\big\{H(t,\bar{x}(t),\bar{u}(t),\psi(t))\big\}\subset D_t^+V(t,\bar{x}(t)),
\end{aligned}
\right.
\end{equation}
where $D_x^{1,-}V,D_x^{1,+}V$ denote the first-order sub- and super-jets of $V$ in the $x$-variable, and $D_t^-V,D_t^+V$ denote the sub- and super-jets of $V$ in the $t$-variable, respectively.

For stochastic optimal control problems, the classical result on the connection between the MP and the DPP (see Bensoussan \cite{Ben82} and Yong and Zhou \cite{YZ99}), is known as $p(t)=-V_x(t,\bar{x}(t)),q(t)=-V_{xx}(t,\bar{x}(t))\sigma(t,\bar{x}(t),\bar{u}(t))$, and $V_t(t,\bar{x}(t))=G(t,\bar{x}(t),\bar{u}(t),-V_x(t,\bar{x}(t)),
-V_{xx}(t,\bar{x}(t)))$ involving an adjoint process pair $(p,q)$ and a generalized Hamiltonian function $G$, where $\sigma$ is the diffusion coefficient. Within the framework of viscosity solution, Zhou \cite{Zhou90-2,Zhou91} showed that
\begin{equation}\label{YZ99}
\left\{
\begin{aligned}
&\{-p(t)\}\times[-P(t),\infty)\subset D_x^{2,+}V(t,\bar{x}(t)),\\
&D_x^{2,-}V(t,\bar{x}(t))\subset\{-p(t)\}\times(-\infty,-P(t)],\\
&[\mathcal{H}(t,\bar{x}(t),\bar{u}(t)),\infty)\subset D_{t+}^{1,+}V(t,\bar{x}(t)),
\end{aligned}
\right.
\end{equation}
where $P$ is the second-order adjoint process (see Peng \cite{Peng90}), the function $\mathcal{H}(t,x,u)$ equals to $G(t,x,u,p(t),P(t))+\mbox{tr}\big\{\sigma(t,x,u)^\top\big[q(t)-P(t)\sigma(t,\bar{x}(t),\bar{u}(t))\big]\big\}$, the notations $D_x^{2,-}V$ (resp. $D_x^{2,+}V$) and $D_{t+}^{1,+}V$ denote the second-order sub- (resp. super-) jets of $V$ in the $x$-variable, and the right super-jet of $V$ in the $t$-variable, respectively.

In this paper, we consider one kind of stochastic recursive optimal control problem, where the cost functional is described by the solution to a {\it backward stochastic differential equation} (BSDE) of the following form
\begin{equation*}
\left\{
\begin{aligned}
-dy(t)&=f(t,y(t),z(t))dt-z(t)dW(t),\ t\in[0,T],\\
  y(T)&=\xi,
\end{aligned}
\right.
\end{equation*}
where the terminal condition $\xi$ is given in advance.
Linear BSDE was introduced by Bismut \cite{Bis78}, to represent the adjoint equation when applying the MP to solve stochastic optimal control problems. Pardoux and Peng \cite{PP90} first studied the adapted solution for the general nonlinear BSDE. Independently, BSDE was involved in Duffie and Epstein \cite{DE92} from economic background, and they presented a stochastic differential formulation of recursive utility which is an extension of the standard additive utility with the instantaneous utility depending not only on the instantaneous consumption rate but also on the future utility. Stochastic recursive optimal control problems have found important applications in mathematical economics, mathematical finance and engineering (see El Karoui, Peng and Quenez \cite{EPQ97}, Wang and Wu \cite{WW09}, etc).

For stochastic recursive optimal control problems, Peng \cite{Peng93} first obtained a local maximum principle when the control domain is convex, by representing the adjoint equation as a {\it forward-backward stochastic differential equation} (FBSDE). Then Xu \cite{Xu95} studied the non-convex control domain case in which the diffusion coefficient does not depends on the control variable. Wu \cite{Wu13} established a general maximum principle for a controlled forward-backward stochastic system where the control domain is non-convex and the diffusion coefficient contains the control variable. The idea of \cite{Wu13} is transferring the original control problem to an equivalent problem with state constraint by understanding the term $z$ appearing in the BSDEs as a control variable. For more general case, the reader is referred to Yong \cite{Yong2010}. Recently, Hu \cite{Hu15} obtained also a general maximum principle for the stochastic recursive optimal control problem, by introducing new and general first- and second-order adjoint equations which are both BSDEs. The result of \cite{Hu15,Wu13,Yong2010}, especially \cite{Hu15}, solves a long-standing open problem of Peng \cite{Peng98} in stochastic control theory, which motivates us to study the connection between this general MP and the DPP. The DPP for stochastic recursive optimal control problems have been solved by Peng \cite{Peng92} (see also Peng \cite{Peng97}), and he proved that the value function is the viscosity solution to a generalized {\it Hamilton-Jacobi-Bellman} (HJB) equation. For further research, the readers are referred to Wu and Yu \cite{WY2008}, Buckdahn and Nie \cite{BN2016} and the references therein.

Concerning the connection between the MP and the DPP for stochastic recursive optimal control problems, Shi \cite{Shi10} first studied (see also Shi and Yu \cite{SY13}) its local form when the control domain is convex and the value function is smooth. The main result is
\begin{equation}\label{Shi10}
\left\{
\begin{aligned}
p^*(t)&=V_x(t,\bar{x}(t))^\top q^*(t),\\
k^*(t)&=\big[V_{xx}(t,\bar{x}(t))\sigma(t,\bar{x}(t),\bar{u}(t))+V_x(t,\bar{x}(t))\times\\
    &\qquad f_z\big(t,\bar{x}(t),-V(t,\bar{x}(t)),-V_x(t,\bar{x}(t))\sigma(t,\bar{x}(t),\bar{u}(t)),\bar{u}(t)\big)\big]q^*(t),
\end{aligned}
\right.
\end{equation}
where $(p^*,q^*,k^*)$ is a triple of adjoint processes and $f$ is the generator of the controlled BSDE. Shi \cite{Shi10} also showed that
\begin{equation}\label{Shi10-2}
\begin{aligned}
  V_t(t,\bar{x}(t))=G\big(t,\bar{x}(t),-V(t,\bar{x}(t),-V_x(t,\bar{x}(t)),-V_{xx}(t,\bar{x}(t)),\bar{u}(t)\big),
\end{aligned}
\end{equation}
where $G$ is the generalized Hamiltonian function of Peng \cite{Peng97}.

However, the above connections (\ref{Shi10}) and (\ref{Shi10-2}) require the smoothness of the value function which does not hold in general, see Example 3.1 of the current paper. This is a major deficiency and important gap of \cite{Shi10,SY13}. This paper will bridge this gap by employing the notions of sub- and super-jets evoked in the viscosity solutions. Moreover, in this paper we study the general case in which the control domain could be non-convex. When the domain of the control is convex, in Nie, Shi and Wu \cite{NSW16} which is the first part of our work we have established the connection between the adjoint process triple $(p^*,q^*,k^*)$ in the MP and the first-order sub- (resp. super-) jets of the value function $V(s,\bar{X}^{t,x;\bar{u}}(s))$, which is
\begin{equation}\label{NieShi16}
\begin{aligned}
D_x^{1,-}V(s,\bar{X}^{t,x;\bar{u}}(s))\subset\{p^*(s)q^*(s)^{-1}\}\subset D_x^{1,+}V(s,\bar{X}^{t,x;\bar{u}}(s)).
\end{aligned}
\end{equation}

As the second part, in the current paper, we study the non-smooth version of the connection between the MP and the DPP for the stochastic recursive optimal control problem with non-convex control domain. We emphasize that in such case we should also study the connection between the second-order sub- and super-jets of the value function and the second-order adjoint process, which is rather difficult. The contributions of this paper are as follows. Firstly, instead of using the adjoint equations of \cite{Peng93} which are FBSDEs, we obtain another form of (\ref{NieShi16}) by using the first-order adjoint equation of \cite{Hu15} which are BSDEs. This is a counterpart of the classical (non-recursive) result of \cite{YZ99}. Secondly, the connection between the second-order adjoint equation of \cite{Hu15} and the second-order sub- (resp. super-) jets of the value function $V$ in the $x$-variable is obtained. See Theorem \ref{main theorem-x} and Theorem \ref{theorem-NSW-ACC2016}. Thirdly, the connection between the right super-differential of the value function $V$ in the $t$-variable and a new $\mathcal{H}_1$-function is derived, where $\mathcal{H}_1$ is although different to the Hamiltonian function $\mathcal{H}$ in the MP of \cite{Hu15}, but by comparing it with the one in the non-recursive case (see \cite{YZ99}), one can see that it is indeed a natural form when extending the classical one of \cite{YZ99} to the recursive case. See Theorem \ref{main theorem-t}. Finally, we also discuss our results from the point view of deriving the MP of \cite{Hu15} directly from the DPP of \cite{Peng92} in smooth case. See Corollary \ref{corollary-DPP to MP}.

Different from the results of \cite{Zhou90-2,Zhou91} for non-recursive stochastic optimal control problems, the recursive structure given by the BSDE arises enormous difficulties for establishing our main results: Theorem 3.1 and Theorem 3.2. In fact, the additional term $Z$ appeared in the BSDE is useful to establish the wellposdness of the BSDE, however it will become a major obstacle when studying the control problems. Since we only have the square estimate for $Z$ from the BSDE theory, i.e. we only have the estimate for $\mathbb{E}\left(\int_0^T|Z_s|^2ds\right)^{p}$, but $|Z_s|^2$ can not be changed to $|Z_s|^{2+\varepsilon}$, for $\varepsilon>0$. This feature makes the classical Taylor expansion arguments in \cite{Zhou90-2,Zhou91,NSW16} fail. We establish our main theorems by introducing the variation to $Y$ and $Z$ motivated by \cite{Hu15}, this is entirely new comparing with \cite{Zhou90-2,Zhou91}. On the other hand, we mention that from the classical results, one may consider that the first order adjoint equation should coincide with the first order derivative of the value function, and the second order adjoint equation should come from the second order derivative of the value function. Then, it should also work in the stochastic recursive optimal control problems, i.e. Theorem 3.1 and Theorem 3.2 holds. From this point of view, one can obtain the adjoint equations which are introduced by \cite{Hu15}. Moreover, conversely, the proof of Theorem 3.1 also gives the illumination for how to construct the variation to $Y$ and $Z$ when proving the general maximum principle. See Section 4 for the details.

Let us mention that it is very important to study the connection between the MP and the DPP for stochastic recursive optimal control problems, because of its valuable applications in mathematical finance and stochastic control theory. For example, on the one hand, in mathematical finance the connection (between the MP and the DPP) can explain well the meaning of the {\it shadow price}, see pp. 254 in \cite{YZ99}. The results obtained in the current paper can help us to understand the shadow price better in the framework of recursive utility. On the other hand, the connection between the MP and the DPP indeed reveals the relationship between the solutions to the adjoint equations and that to the HJB equation. As well known, the adjoint equations are usually {\it fully-coupled} FBSDEs which are very hard to solve in general. The relationship between the adjoint processes and the solution to the HJB equation give us a heuristic method, to obtain the solutions to the adjoint equations by the corresponding HJB equation which can be solved analytically or numerically. Moreover, as shown in Corollary \ref{corollary-DPP to MP}, such connection also allows us to obtain the MP (which is difficult for the general case and is just solved recently) from the DPP (which have been solved since 1992 by Peng \cite{Peng92}). These relations are also helpful to find the candidate optimal controls (see Remark 3.1) and to establish the stochastic verification theorem. All the above important issues motivate us to study the connection between the MP and the DPP for the general stochastic recursive optimal control problems in this paper.

The rest of this paper is organized as follows. In Section 2, we state our problem and give some preliminary results about the MP and the DPP. Section 3 exhibits the main results of this paper, specially we give the connections between the value function and the adjoint processes within the framework of viscosity solution; we derive the general first and second order adjoint equations; we prove the general MP from the DPP by using the results of this paper. We also give several interesting examples to explain the related results. Section 4 devotes to the proof of the main results. Finally, in Section 5, we give the concluding remarks.

In this paper, $\mathbf{R}^n$ denotes the $n$-dimensional Euclidean space with scalar product $\langle\cdot,\cdot\rangle$ and norm $|\cdot|$, $\mathcal{S}^n$ denotes the $n\times n$ symmetric matrix space, $Df,D^2f$ denotes the gradient and the Hessian matrix of the differentiable function $f$ respectively, $\top$ appearing as superscript denotes the transpose of a matrix, and $C>0$ denotes a generic constant which may take different values in different places.

\section{Problem Statement and Preliminaries}

Let $T>0$ be finite and $\mathbf{U}\subset\mathbf{R}^k$ be nonempty. Given $t\in[0,T)$, we denote $\mathcal{U}^w[t,T]$ the set of all 5-tuples $(\Omega,\mathcal{F},\mathbf{P},W(\cdot);\\u(\cdot))$ satisfying the following conditions: (i) $(\Omega,\mathcal{F},\mathbf{P})$ is a complete probability space; (ii) $\{W(s)\}_{s\geq t}$ is a one-dimensional standard Brownian motion defined on $(\Omega,\mathcal{F},\mathbf{P})$ over $[t,T]$ (with $W(t)=0$, a.s.), and $\mathcal{F}_s^t=\sigma\{W(r);t\leq r\leq s\}$ augmented by all the $\mathbf{P}$-null sets in $\mathcal{F}$; (iii) $u:[t,T]\times\Omega\rightarrow\mathbf{U}$ is an $\{\mathcal{F}_s^t\}_{s\geq t}$-adapted process on $(\Omega,\mathcal{F},\mathbf{P})$.

We write $(\Omega,\mathcal{F},\mathbf{P},W(\cdot);u(\cdot))\in\mathcal{U}^w[t,T]$, but occasionally we will write only $u(\cdot)\in\mathcal{U}^w[t,T]$ if no ambiguity exists. We denote by $L^2(\Omega,\mathcal{F}^t_T;\mathbf{R}^n)$ the space of all $\mathbf{R}^n$-valued $\mathcal{F}^t_T$-measurable random variables
$\xi$ such that $\mathbb{E}[|\xi|^2]<\infty$, $S^2_\mathcal{F}([t,T];\mathbf{R}^n)$ denotes the space of all $\mathbf{R}^n$-valued $\mathcal{F}^t_s$-adapted process $\phi(\cdot)$ such that $\mathbb{E}\big[\sup_{t\leq s\leq T}|\phi(s)|^2\big]\\<\infty$, and $L^2_\mathcal{F}([t,T];\mathbf{R}^n)$ denotes the space of all $\mathbf{R}^n$-valued $\mathcal{F}^t_s$-adapted processes $\psi(\cdot)$ such that $\mathbb{E}\big[\int_t^T|\psi(s)|^2ds\big]<\infty$.

For any $(t,x)\in[0,T)\times\mathbf{R}^n$, we consider the state $X^{t,x;u}(\cdot)\in{\mathbf{R}}^n$ given by the following controlled {\it stochastic differential equation} (SDE):
\begin{equation}\label{controlled SDE}
\left\{
\begin{aligned}
 dX^{t,x;u}(s)&=b(s,X^{t,x;u}(s),u(s))ds+\sigma(s,X^{t,x;u}(s),u(s))dW(s),\ s\in[t,T],\\
  X^{t,x;u}(t)&=x.
\end{aligned}
\right.
\end{equation}
Here $b:[0,T]\times\mathbf{R}^n\times\mathbf{U}\rightarrow\mathbf{R}^n,\sigma:[0,T]\times\mathbf{R}^n\times\mathbf{U}\rightarrow\mathbf{R}^n$ are given functions satisfying

\vspace{1mm}

\noindent({\bf H1})\ $b$ and $\sigma$ are uniformly continuous in $(s,x,u)$, Lipschitz continuous and linear growth in $x$.

For any $(t,x)\in[0,T)\times\mathbf{R}^n$ and given $u(\cdot)\in\mathcal{U}^w[t,T]$, under ({\bf H1}), SDE (\ref{controlled SDE}) has a unique solution $X^{t,x;u}(\cdot)\in S^2_\mathcal{F}([t,T];\mathbf{R}^n)$ (see \cite{YZ99}). We refer to such $u(\cdot)\in\mathcal{U}^w[t,T]$ as an admissible control and $(X^{t,x;u}(\cdot),u(\cdot))$ as an admissible pair.

Next, we introduce the following controlled BSDE coupled with (\ref{controlled SDE}):
\begin{equation}\label{controlled BSDE}
\left\{
\begin{aligned}
-dY^{t,x;u}(s)&=f(s,X^{t,x;u}(s),Y^{t,x;u}(s),Z^{t,x;u}(s),u(s))ds\\
              &\quad -Z^{t,x;u}(s)dW(s),\ s\in[t,T],\\
  Y^{t,x;u}(T)&=\phi(X^{t,x;u}(T)),
\end{aligned}
\right.
\end{equation}
where $f:[0,T]\times\mathbf{R}^n\times\mathbf{R}\times\mathbf{R}\times\mathbf{U}\rightarrow\mathbf{R}$ and $\phi:\mathbf{R}^n\rightarrow\mathbf{R}$ are given functions satisfying

\vspace{1mm}

\noindent({\bf H2})\ $f$ and $\phi$ are uniformly continuous in $(s,x,y,z,u)$, Lipschitz continuous in $(x,y,z)$ and linear growth in $x$.

\vspace{1mm}

\begin{lemma}\label{lemma-PP1990}
(\cite{PP90,Peng97}) Let {\bf (H1), (H2)} hold. For given $(t,x)\in[0,T)\times\mathbf{R}^n$, $u(\cdot)\in\mathcal{U}^w[t,T]$, and $X^{t,x;u}(\cdot)\in S^2_\mathcal{F}([t,T];\mathbf{R}^n)$ is the unique solution to (\ref{controlled SDE}). Then BSDE (\ref{controlled BSDE}) admits a unique solution $(Y^{t,x;u}(\cdot),Z^{t,x;u}(\cdot))\in S^2_\mathcal{F}([t,T];\mathbf{R})\times L^2_\mathcal{F}([t,T];\mathbf{R})$. Moreover, the following estimate holds:
\begin{equation*}\label{basic estimate of BSDE}
{\small\begin{aligned}
\mathbb{E}\Big[\sup\limits_{t\leq s\leq T}\big|Y^{t,x;u}(s)\big|^2\Big|\mathcal{F}_s^t\Big]
\leq C\mathbb{E}\Big[\big|\phi(X^{t,x;u}(T))\big|^2+\int_t^T\big|f\big(s,X^{t,x;u}(s),0,0,u(s)\big)\big|^2ds\Big|\mathcal{F}_s^t\Big].
\end{aligned}}
\end{equation*}
\end{lemma}

{\bf Remark 2.1}\ The wellposedness of BSDE can be established in the case that $f$ is quadratic growth in $Z$ and with bounded or unbounded terminal condition, see Kobylanski \cite{K00}, Briand and Hu \cite{BH06,BH08} and the reference therein. For stochastic non-recursive optimal control problems related to BSDEs with quadratic growth, the reader is referred to Fuhrman, Hu and Tessitore \cite{FHT06}, for example. However, to our best knowledge, until now there are no works on stochastic recursive control problems whose cost functional is given by quadratic BSDEs. Concerning the assumptions of the coefficients on $x$, we mention that polynomial growth of $f$ and $\phi$ on $x$ is enough to obtain the wellposedness of BSDE, and continuity assumption of $f$ and $\phi$ on $x$ is sufficient to show that the value function is a viscosity solution to HJB equation. When proving the uniqueness of the viscosity solution to HJB equation, usually $f$ and $\phi$ are assumed to be Lipschitz on $x$. The readers are referred to Wu and Yu \cite{WY2008}, Buckdahn and Nie \cite{BN2016} and the references therein. In our paper, in order to focus on our problems at hand  and also to make the notations easier, we prefer to keep using the standard Lipschitz assumption for the coefficients on $(x,y,z)$.

Given $u(\cdot)\in\mathcal{U}^w[t,T]$, we introduce the cost functional as
\begin{equation}\label{cost functional}
J(t,x;u(\cdot)):=-Y^{t,x;u}(s)|_{s=t},\quad(t,x)\in[0,T)\times\mathbf{R}^n.
\end{equation}
and our stochastic recursive optimal control problem is the following.

\vspace{1mm}

{\bf Problem (SROCP).}\quad For given $(t,x)\in[0,T)\times\mathbf{R}^n$, to minimize (\ref{cost
functional}) subject to (\ref{controlled SDE})$\sim$(\ref{controlled BSDE}) over $\mathcal{U}^w[t,T]$.

We define the value function
\begin{equation}\label{value function}
\left\{
\begin{aligned}
&V(t,x):=\inf\limits_{u(\cdot)\in\mathcal{U}^w[t,T]}J(t,x;u(\cdot)),\ (t,x)\in[0,T)\times\mathbf{R}^n,\\
&V(T,x)=-\phi(x),\quad x\in\mathbf{R}^n.
\end{aligned}
\right.
\end{equation}
Any $\bar{u}(\cdot)\in\mathcal{U}^w[t,T]$ achieves the above infimum is called an optimal control, and the corresponding solution triple $(\bar{X}^{t,x;\bar{u}}(\cdot),\bar{Y}^{t,x;\bar{u}}(\cdot),\bar{Z}^{t,x;\bar{u}}(\cdot))$ is called an optimal trajectory. We refer to $(\bar{X}^{t,x;\bar{u}}(\cdot),\bar{Y}^{t,x;\bar{u}}(\cdot),\bar{Z}^{t,x;\bar{u}}(\cdot),\bar{u}(\cdot))$ as an optimal quadruple.

\vspace{1mm}

{\bf Remark 2.2}\ Since the functions $b,\sigma,f$ and $g$ are all deterministic, from Proposition 5.1 of \cite{Peng97}, we know that under ({\bf H1}), ({\bf H2}), the cost functional (\ref{cost functional}) is deterministic. Thus our definition (\ref{value function}) is meaningful.

We introduce the following generalized HJB equation:
\begin{equation}\label{HJB equation}
\left\{
\begin{aligned}
 &-v_t(t,x)+\sup\limits_{u\in\mathbf{U}}G(t,x,-v(t,x),-v_x(t,x),-v_{xx}(t,x),u)=0,\ (t,x)\in[0,T)\times\mathbf{R}^n,\\
 &v(T,x)=-\phi(x),\quad \forall x\in\mathbf{R}^n,
\end{aligned}
\right.
\end{equation}
where the generalized Hamiltonian function $G:[0,T]\times\mathbf{R}^n\times\mathbf{R}\times\mathbf{R}^n\times\mathcal{S}^n\times\mathbf{U}\rightarrow\mathbf{R}$ is defined as
\begin{equation}\label{generalized Hamiltonian}
{\small\begin{aligned}
 G(t,x,r,p,A,u):=\frac{1}{2}\langle A\sigma(t,x,u),\sigma(t,x,u)\rangle+\langle p,b(t,x,u)\rangle+f(t,x,r,\sigma(t,x,u)^\top p,u).
\end{aligned}}
\end{equation}

We introduce the following definition of the viscosity solution to HJB equation (\ref{HJB equation}), which can be found in Crandall, Ishii and Lions \cite{CIL92}.

\vspace{1mm}

\begin{definition}\label{definition-viscosity solution}
(\cite{CIL92}) (i) A function $v\in C([0,T]\times\mathbb{R}^n)$ is called a viscosity sub- (resp. super-) solution to (\ref{HJB equation}) if
$v(T,x)\leq -\phi(x)$ (resp. $v(T,x)\geq -\phi(x)$), for all $x\in\mathbb{R}^n$,
and for any $\varphi\in C^{1,2}([0,T]\times\mathbb{R}^n)$ such that $v-\varphi$ attains a local maximum (resp. minmum) at $(t,x)\in[0,T)\times\mathbb{R}^n$, we have
\begin{equation*}
\begin{aligned}
&-\varphi_t(t,x)+\sup\limits_{u\in\mathbf{U}}G\big(t,x,-v(t,x),-\varphi_x(t,x),-\varphi_{xx}(t,x),u\big)\leq0\\
&(resp.\ -\varphi_t(t,x)+\sup\limits_{u\in\mathbf{U}}G\big(t,x,-v(t,x),-\varphi_x(t,x),-\varphi_{xx}(t,x),u\big)\geq0).
\end{aligned}
\end{equation*}
(ii) A function $v\in C([0,T]\times\mathbb{R}^n)$ is called a viscosity solution to (\ref{HJB equation}) if it is both a viscosity subsolution and viscosity supersolution to (\ref{HJB equation}).
\end{definition}

\vspace{1mm}

{\bf Remark 2.3}\ The viscosity solution to HJB equation (\ref{HJB equation}) can be equivalently defined in the language of sub- and
super-jets, see \cite{CIL92, YZ99}.

\vspace{1mm}

\begin{proposition}\label{proposition-Peng1997}
(\cite{Peng97}) Let {\bf (H1), (H2)} hold. Let $V(\cdot,\cdot)$ be defined by (\ref{value function}). Then for any $t\in[0,T]$ and $x,x'\in\mathbf{R}^n$, we have
\begin{equation}\label{value function: regularity}
\begin{aligned}
\mbox{(i)}\ |V(t,x)-V(t,x')|\leq C|x-x'|,\quad\mbox{(ii)}\ |V(t,x)|\leq C(1+|x|).
\end{aligned}
\end{equation}
Moreover, $V(\cdot,\cdot)$ is the unique viscosity solution to (\ref{HJB equation}).
\end{proposition}

To conveniently state the maximum principle,  for given $(t,x)\in[0,T)\times\mathbf{R}^n$, we rewrite the equations (\ref{controlled SDE}) and (\ref{controlled BSDE}) as the following controlled FBSDE:
\begin{equation}\label{controlled FBSDE}
\left\{
{\small\begin{aligned}
 dX^{t,x;u}(s)&=b(s,X^{t,x;u}(s),u(s))ds+\sigma(s,X^{t,x;u}(s),u(s))dW(s),\\
-dY^{t,x;u}(s)&=f(s,X^{t,x;u}(s),Y^{t,x;u}(s),Z^{t,x;u}(s)u(s))ds-Z^{t,x;u}(s)dW(s),\ s\in[t,T],\\
  X^{t,x;u}(t)&=x,\ Y^{t,x;u}(T)=\phi(X^{t,x;u}(T)).
\end{aligned}}
\right.
\end{equation}
Moreover, we introduce the following assumption as in \cite{Hu15}. Recall that $Df,D^2f$ denotes, respectively, the gradient and the Hessian matrix of the differentiable function $f$ with respect to the variable $(x,y,z)$.

\vspace{1mm}

\noindent({\bf H3})\ The functions $b,\sigma,\phi$ and $f$ are twice continuously differentiable in $(x,y,z)$, the derivatives  $b_x,b_{xx},\sigma_x,\sigma_{xx},\phi_x,\phi_{xx},Df$ and $D^2f$ are continuous in $(s,x,y,z,u)$ and uniformly bounded. Moreover, there exists a modulus of continuity $\overline{\omega}:[0,\infty]\rightarrow[0,\infty]$ for $\varphi=b_{xx},\sigma_{xx},\phi_{xx}$ and $\psi=D^2f$, such that for arbitrary $s\in[0,T]$, $u\in \mathbf{U}$, $x_1,x_2\in\mathbf{R}^n$, $y_1,y_2,z_1,z_2\in\mathbf{R}$,
\begin{equation*}
\begin{aligned}
&|\varphi(s,x_1,u)-\varphi(s,x_2,u)|\leq \overline{\omega}(|x_1-x_2|),\\
&|\psi(s,x_1,y_1,z_1,u)-\psi(s,x_2,y_2,z_2,u)|\leq \overline{\omega}(|x_1-x_2|+|y_1-y_2|+|z_1-z_2|).
\end{aligned}
\end{equation*}

Let $(\bar{X}^{t,x;\bar{u}}(\cdot),\bar{Y}^{t,x;\bar{u}}(\cdot),\bar{Z}^{t,x;\bar{u}}(\cdot),\bar{u}(\cdot))$ be an optimal quadruple. For given $(t,x)\in[0,T)\times\mathbf{R}^n$,  we introduce the following first- and second-order adjoint equation (\cite{Hu15})
\begin{equation}\label{adjoint equation-first order}
\left\{
{\small\begin{aligned}
-dp(s)&=\Big\{\bar{f}_y(s)p(s)+\big[\bar{f}_z(s)\bar{\sigma}_x(s)^\top+\bar{b}_x(s)^\top\big]p(s)+\bar{f}_z(s)q(s)\\
      &\quad\ +\bar{\sigma}_x(s)^\top q(s)+\bar{f}_x(s)\Big\}ds-q(s)dW(s),\ s\in[t,T],\\
  p(T)&=\phi_x(\bar{X}^{t,x;\bar{u}}(T)),
\end{aligned}}
\right.
\end{equation}
\begin{equation}\label{adjoint equation-second order}
\left\{
{\small\begin{aligned}
-dP(s)&=\Big\{\bar{f}_y(s)P(s)+\big[\bar{f}_z(s)\bar{\sigma}_x(s)^\top+\bar{b}_x(s)^\top\big]P(s)+P(s)\big[\bar{f}_z(s)\bar{\sigma}_x(s)+\bar{b}_x(s)\big]\\
      &\quad\ +\bar{\sigma}_x(s)^\top P(s)\bar{\sigma}_x(s)+\bar{f}_z(s)Q(s)+\bar{\sigma}_x(s)^\top Q(s)+Q(s)\bar{\sigma}_x(s)+\bar{b}_{xx}(s)^\top p(s)\\
      &\quad\ +\bar{\sigma}_{xx}(s)^\top\big[\bar{f}_z(s)p(s)+q(s)\big]+\big[I_{n\times n},p(s),\bar{\sigma}_x(s)^\top p(s)+q(s)\big]D^2\bar{f}(s)\\
      &\quad\ \cdot\big[I_{n\times n},p(s),\bar{\sigma}_x(s)^\top p(s)+q(s)\big]^\top\Big\}ds-Q(s)dW(s),\ s\in[t,T],\\
  P(T)&=\phi_{xx}(\bar{X}^{t,x;\bar{u}}(T)),
\end{aligned}}
\right.
\end{equation}
respectively. Here $\bar{b}(s):=b(s,\bar{X}^{t,x;\bar{u}}(s),\bar{u}(s)),\bar{\sigma}(s):=\sigma(s,\bar{X}^{t,x;\bar{u}}(s),\bar{u}(s)),\bar{f}(s):=f(s,\bar{X}^{t,x;\bar{u}}(s),\bar{Y}^{t,x;\bar{u}}(s),\bar{Z}^{t,x;\bar{u}}(s),\bar{u}(s))$, and similar notations are used for all their derivatives, and $p(s)=(p^1(s),p^2(s),\cdots,p^n(s))^\top,\bar{b}_{xx}(s)^\top p(s)=\sum_{i=1}^{n}p^i(s)b^i_{xx}(s)$, and for $\bar{\sigma}_{xx}(s)^\top\big[\bar{f}_z(s)p(s)+q(s)\big]$ similarly. The term $\big[I_{n\times n},p(s),\bar{\sigma}_x(s)^\top p(s)+q(s)\big]$ denotes the $n\times (n+2)$-matrix combining $I_{n\times n}$, $p(s)$ and $\bar{\sigma}_x(s)^\top p(s)+q(s)$. Noticing that $D^2\bar{f}$ is a $(n+2)\times (n+2)$-matrix, thus the term $\big[I_{n\times n},p(s),\bar{\sigma}_x(s)^\top p(s)+q(s)\big]D^2\bar{f}(s)\big[I_{n\times n},p(s),\bar{\sigma}_x(s)^\top p(s)+q(s)\big]^\top$ is well defined.

By {\bf (H3)}, it is easy to verify that the above BSDEs admit the unique solution $(p(\cdot),q(\cdot))\in S^2_\mathcal{F}([t,T];\mathbf{R}^n)\times L^2_\mathcal{F}([t,T];\mathbf{R}^n)$ and $(P(\cdot),Q(\cdot))\in S^2_\mathcal{F}([t,T];\mathcal{S}^n)\times L^2_\mathcal{F}([t,T];\mathcal{S}^n)$, respectively.

\vspace{1mm}

\begin{proposition}\label{proposition-Hu2015}
(\cite{Hu15}) Let {\bf (H1), (H2), (H3)} hold
and $(t,x)\in[0,T)\times\mathbf{R}^n$ be fixed. Suppose that $\bar{u}(\cdot)$ is an optimal control for {\bf Problem (SROCP)},
$(\bar{X}^{t,x;\bar{u}}(\cdot),\bar{Y}^{t,x;\bar{u}}(\cdot),\bar{Z}^{t,x;\bar{u}}(\cdot))$
is the optimal trajectory, and $(p(\cdot),q(\cdot)),(P(\cdot),Q(\cdot))$ satisfy (\ref{adjoint equation-first order}) and (\ref{adjoint equation-second order}), respectively. Then
\begin{equation}\label{maximum condition}
\begin{aligned}
&\mathcal{H}(s,\bar{X}^{t,x;\bar{u}}(s),\bar{Y}^{t,x;\bar{u}}(s),\bar{Z}^{t,x;\bar{u}}(s),\bar{u}(s),p(s),q(s),P(s))\\
=&\max\limits_{u\in\mathbf{U}}\mathcal{H}(s,\bar{X}^{t,x;\bar{u}}(s),\bar{Y}^{t,x;\bar{u}}(s),\bar{Z}^{t,x;\bar{u}}(s),u,p(s),q(s),P(s)),\
\end{aligned}
\end{equation}
a.e. $s\in[t,T],\mathbf{P}\mbox{-}$a.s., where the Hamiltonian function $\mathcal{H}:[0,T]\times\mathbf{R}^n\times\mathbf{R}\times\mathbf{R}\times\mathbf{U}\times\mathbf{R}^n\times\mathbf{R}^n\times\mathcal{S}^n
\rightarrow\mathbf{R}$ is defined as
\begin{equation}\label{Hamiltonian}
\begin{aligned}
&\mathcal{H}(t,x,y,z,u,p,q,P):=f\big(t,x,y,z+\langle p,\sigma(t,x,u)-\sigma(t,\bar{x},\bar{u})\rangle,u\big)+\langle p,b(t,x,u)\rangle\\
&\quad+\langle q,\sigma(t,x,u)\rangle+\frac{1}{2}\big\langle P(\sigma(t,x,u)-\sigma(t,\bar{x},\bar{u})),\sigma(t,x,u)-\sigma(t,\bar{x},\bar{u})\big\rangle.
\end{aligned}
\end{equation}
\end{proposition}

{\bf Remark 2.4}\ Notice that Proposition \ref{proposition-Hu2015} is proved by \cite{Hu15} in its strong formulation. However, as pointed out in \cite{YZ99}, since for the maximum principle only necessary conditions of optimality are concerned, an optimal quadruple is given as a starting point (no matter whether in the strong or weak formulation, here we mention that the strong/weak formulation concerns on the structure of the control variable, not on the strong/weak solution of SDEs, for more details, see Chapter 2, Section 4 in \cite{YZ99}), and all the results are valid for this given optimal quadruple on the probability space it attached to.

{\bf Remark 2.5}\ In fact, the first-order adjoint equation (\ref{adjoint equation-first order}) can be derived by comparing the result in \cite{YZ99} with the result in \cite{NSW16} which using the adjoint FBSDE of \cite{Peng93}, for more details see Theorem \ref{theorem-NSW-ACC2016} and Remark 3.3. Moreover, under some conditions we can also obtain the second-order adjoint equation (\ref{adjoint equation-second order}) by comparing the results between the recursive case and the non-recursive case, see Corollary \ref{corollary-smooth case}. Finally, using the relationship established in the next section, we can prove Proposition \ref{proposition-Hu2015} directly from the DPP of \cite{Peng92}, see Corollary \ref{corollary-DPP to MP}.

\section{Main Results}

\subsection{Differentials in the spatial variable}

In this subsection, we extend (\ref{Shi10}) to the case with non-convex control domain and non-smooth value function.

Let us recall the notions of the second-order super- and sub-jets in the spatial variable $x$,  see \cite{CIL92,YZ99}. For $v\in C([0,T]\times\mathbb{R}^n)$ and $(t,\hat{x})\in[0,T]\times\mathbb{R}^n$, we define
\begin{equation}\label{second-order super- and sub-jets}
\left\{
\begin{aligned}
D_x^{2,+}v(t,\hat{x})&:=\Big\{(p,P)\in\times\mathbf{R}^n\times\mathcal{S}^n\big|v(t,x)\leq v(t,\hat{x})+\langle p,x-\hat{x}\rangle\\
                     &\qquad+\frac{1}{2}(x-\hat{x})^\top P(x-\hat{x})+o(|x-\hat{x}|^2),\mbox{ as }x\rightarrow\hat{x}\Big\},\\
D_x^{2,-}v(t,\hat{x})&:=\Big\{(p,P)\in\times\mathbf{R}^n\times\mathcal{S}^n\big|v(t,x)\geq v(t,\hat{x})+\langle p,x-\hat{x}\rangle\\
                     &\qquad+\frac{1}{2}(x-\hat{x})^\top P(x-\hat{x})+o(|x-\hat{x}|^2),\mbox{ as }x\rightarrow\hat{x}\Big\}.
\end{aligned}
\right.
\end{equation}

\vspace{1mm}

\begin{theorem}\label{main theorem-x}
Let {\bf (H1), (H2), (H3)} hold, $(t,x)\in[0,T)\times\mathbf{R}^n$ be fixed. Suppose that $\bar{u}(\cdot)$ is an optimal control for {\bf Problem (RSOCP)}, $V(\cdot,\cdot)$ is the value function, and
$(\bar{X}^{t,x;\bar{u}}(\cdot),\bar{Y}^{t,x;\bar{u}}(\cdot),\bar{Z}^{t,x;\bar{u}}(\cdot))$
is the optimal trajectory. Let $(p(\cdot),q(\cdot))\in S^2_\mathcal{F}([t,T];\mathbf{R}^n)\times L^2_\mathcal{F}([t,T];\mathbf{R}^n)$, $(P(\cdot),Q(\cdot))\in S^2_\mathcal{F}([t,T];\mathcal{S}^n)\times L^2_\mathcal{F}([t,T];\mathcal{S}^n)$ satisfy the first- and second-order adjoint equations (\ref{adjoint equation-first order}), (\ref{adjoint equation-second order}), respectively. Then
\begin{equation}\label{connection-second order-x}
\begin{aligned}
&\{-p(s)\}\times[-P(s),\infty)\subseteq D_x^{2,+}V(s,\bar{X}^{t,x;\bar{u}}(s)),\\
&D_x^{2,-}V(s,\bar{X}^{t,x;\bar{u}}(s))\subseteq\{-p(s)\}\times(-\infty,-P(s)],\ \text{ for all } s\in[t,T],\mathbf{P}\mbox{-}a.s.
\end{aligned}
\end{equation}
\end{theorem}

We postpone the proof of Theorem 3.1 to the next section. Let us first give some remarks, examples and corollaries.

\vspace{1mm}

{\bf Remark 3.1}\ Theorem \ref{main theorem-x} also provides a possible method to find the candidate optimal controls. With the stochastic verification theorem of Zhang \cite{Zhang12}, we can solve the original {\bf Problem (RSOCP)} as the following steps. Step 1, solve the corresponding HJB equation (\ref{HJB equation}) to find its unique viscosity solution. Step 2, choose a $\bar{u}(\cdot)\in\mathcal{U}^w[t,T]$ such that (\ref{connection-second order-x}) holds, which is a candidate optimal control. We first calculate the lower bound and super bound of first-component of $D_x^{2,+}V(s,\bar{X}^{t,x;\bar{u}}(s))$ and the lower bound of the second-component of $D_x^{2,+}V(s,\bar{X}^{t,x;\bar{u}}(s))$. Then using the comparison theorem of BSDEs, by comparing adjoint equations (\ref{adjoint equation-first order}) and (\ref{adjoint equation-second order}) with such bounds respectively, we can find some candidate optimal controls with the help of (\ref{connection-second order-x}). Step 3, apply the stochastic verification theorem  to verify the optimal control. The following is an example to show the above solving process.

\vspace{1mm}

{\bf Example 3.1}\quad Let us now consider the following controlled SDE ($n=1$):
\begin{equation}\label{controlled SDE example 1}
\left\{
\begin{aligned}
 dX^{t,x;u}(s)&=\big[X^{t,x;u}(s)+X^{t,x;u}(s)u(s)\big]ds+X^{t,x;u}(s)u(s)dW(s),\ s\in[t,T],\\
  X^{t,x;u}(t)&=x,
\end{aligned}
\right.
\end{equation}
with $\mathbf{U}=[-1,0]\bigcup[1,2]$. The cost functional is defined by (\ref{cost functional}) with
\begin{equation}\label{controlled BSDE example 1}
\left\{
\begin{aligned}
-dY^{t,x;u}(s)&=-Z^{t,x;u}(s)u(s)ds-Z^{t,x;u}(s)dW(s),\ s\in[t,T],\\
  Y^{t,x;u}(T)&=X^{t,x;u}(T).
\end{aligned}
\right.
\end{equation}
The corresponding generalized HJB equation reads as
\begin{equation}\label{HJB equation example 1}
\left\{
\begin{aligned}
 &\sup\limits_{u\in\mathbf{U}}\Big\{-\frac{1}{2}v_{xx}(t,x)x^2u^2+u^2v_x(t,x)x-uv_x(t,x)x\Big\}\\
 &\qquad-v_t(t,x)-v_x(t,x)x=0,\qquad(t,x)\in[0,T)\times\mathbf{R}^n,\\
 &v(T,x)=-x,\quad \text{ for all } x\in\mathbf{R}^n.
\end{aligned}
\right.
\end{equation}
It is easy to verify that, the viscosity solution to (\ref{HJB equation example 1}) is given by
\begin{equation}\label{value function-exapmle 1}
\begin{aligned}
V(t,x)=\left\{\begin{array}{lc}-e^{t-T}x,&\mbox{if }x\leq0,\\
                               -e^{T-t}x,&\mbox{if }x>0,\end{array}\right.
\end{aligned}
\end{equation}
which obviously satisfying (\ref{value function: regularity}). Thus, by the uniqueness of the viscosity solution, $V$ coincides with the value function. Moreover, the first- and second order adjoint equation are
\begin{equation}\label{first-order adjoint equation example 1}
\left\{
\begin{aligned}
-dp(s)&=\big[-|\bar{u}(s)|^2p(s)+p(s)+\bar{u}(s)p(s)\big]ds-q(s)dW(s),\\
  p(T)&=1,
\end{aligned}
\right.
\end{equation}
\begin{equation}\label{second-order adjoint equation example 1}
\left\{
\begin{aligned}
-dP(s)&=\big\{[-|\bar{u}(s)|^2+2+2\bar{u}(s)]P(s)-\bar{u}(s)Q(s)\big\}ds-Q(s)dW(s),\\
  P(T)&=0.
\end{aligned}
\right.
\end{equation}
respectively. Note that $(P(s),Q(s))=0$, for any $\bar{u}(\cdot)\in\mathcal{U}^{w}[t,T]$.

Let us first consider the case that $x\neq 0$. Suppose that the optimal control is $\bar{u}(\cdot)$. If $x\leq0$, from the comparison theorem of SDEs, we have $\bar{X}^{t,x;\bar{u}}(s)\leq0$, for all $s\in[t,T]$. Then $V(s,\bar{X}^{t,x;\bar{u}}(s))=-e^{s-T}\bar{X}^{t,x;\bar{u}}(s)$, which yields that
\begin{equation*}
D_x^{2,+}V(s,\bar{X}^{t,x;\bar{u}}(s))=\{-e^{s-T}\}\times[0,\infty), \quad \text{for all } s\in[t,T].
\end{equation*}
From relation (\ref{connection-second order-x}), we should have $p(s)=e^{s-T}$, for all $s\in[t,T]$. Consequently $q(s)=0$ and $\bar{u}(s)=-1$ or $2$ for all $s\in[t,T]$. This tell us that the candidate optimal control is $\bar{u}(s)=-1$ or $2$. Conversely, using the stochastic verification theorem in \cite{SY13} (which works also for the non-convex control domain since their proof mainly depends on Proposition \ref{proposition-Peng1997}), one can find the optimal control is indeed $\bar{u}(s)=-1$ or $2$, then the unique solution to (\ref{first-order adjoint equation example 1}) is $(p(s),q(s))=(e^{s-T},0)$, for $s\in[t,T]$. Thus (\ref{connection-second order-x}) holds.

If $x\ge0$, from the comparison theorem of SDEs, we have $\bar{X}^{t,x;\bar{u}}(s)\ge0$, for all $s\in[t,T]$. Then $V(s,\bar{X}^{t,x;\bar{u}}(s))=-e^{T-s}\bar{X}^{t,x;\bar{u}}(s)$, which yields that
\begin{equation*}
D_x^{2,+}V(s,\bar{X}^{t,x;\bar{u}}(s))=\{-e^{T-s}\}\times[0,\infty), \quad \text{for all } s\in[t,T].
\end{equation*}
From relation (\ref{connection-second order-x}), we should have $p(s)=e^{T-s}$ for all $s\in[t,T]$, which means that $p(s)=e^{T-s}$ solves uniquely (\ref{first-order adjoint equation example 1}). Consequently $q(s)=0$, $\bar{u}(s)=0$ or $1$ for $s\in[t,T]$. This tell us that the candidate optimal control is $\bar{u}(s)=0$ or $1$. Similarly as above, using the stochastic verification theorem, one can find the optimal control is $\bar{u}(s)=0$ or $1$. Then the unique solution to (\ref{first-order adjoint equation example 1}) is $(p(s),q(s))=(e^{T-s},0)$ and (\ref{connection-second order-x}) holds.

Now let us consider that the initial state $x=0$. We have $\bar{X}^{t,x;\bar{u}}(\cdot)=0$ for any $u\in \mathcal{U}^{w}[t,T]$. From (\ref{value function-exapmle 1}), we have
\begin{equation*}
D_x^{2,-}V(s,\bar{X}^{t,x;\bar{u}}(s))=\emptyset,\quad D_x^{2,+}V(s,\bar{X}^{t,x;\bar{u}}(s))=[-e^{T-s},-e^{s-T}]\times[0,\infty).
\end{equation*}
An optimal control $\bar{u}(\cdot)$ should keep that (\ref{connection-second order-x}) holds where $(p(s),q(s))$ and $(P(s),Q(s))$ solves respectively (\ref{first-order adjoint equation example 1}) and (\ref{second-order adjoint equation example 1}). Thus we should have $e^{s-T}=:p_1(s)\leq p(s)\leq p_2(s):=e^{T-s}$, and $P(s)\leq P_1(s):=0$. On the one hand, since $P(s)=0$ for any $\bar{u}(\cdot)$, thus any $\bar{u}(\cdot)$ can make sure that $P(s)\leq P_1(s)$. On the other hand, it is obvious that $-dp_1(s)=-p_1(s)ds-0dW(s)$, $p_1(T)=1$, then comparing this with (\ref{first-order adjoint equation example 1}) and using the comparison theorem of BSDEs, to make sure $p_1(s)\leq p(s)$, it is sufficient to have $-(\bar{u}(s))^2+1+\bar{u}(s)\ge -1$ which always holds for $u(\cdot)\in\mathbf{U}=[-1,0]\bigcup[1,2]$. Similarly, $-dp_2(s)=p_2(s)ds-0dW(s)$, $p_2(T)=1$, then comparing this with (\ref{first-order adjoint equation example 1}) and using the comparison theorem of BSDEs, to make sure $p(s)\leq p_2(s)$, it is sufficient to have $-(\bar{u}(s))^2+1+\bar{u}(s)\leq 1$ which always holds for $u(\cdot)\in\mathbf{U}=[-1,0]\bigcup[1,2]$. Therefore, any $\bar{u}(\cdot)\in \mathcal{U}^{w}[t,T]$ is a candidate optimal control.

Conversely, it can be verified that $(\bar{X}^{t,x;\bar{u}}(\cdot),\bar{u}(\cdot))=(0,\bar{u}(\cdot))$ is an optimal pair (one can check easily that when $x=0$, the cost functional is always zero no mater what the control is). Moreover, using comparison theorem of BSDEs and $\mathbf{U}=[-1,0]\bigcup[1,2]$, the same as above we can show that $e^{s-T}\leq p(s)\leq e^{T-s}$, for any $\bar{u}(\cdot)\in\mathcal{U}^{w}[t,T]$. Noticing that $P(s)=0$ for any $\bar{u}(\cdot)\in\mathcal{U}^{w}[t,T]$. Thus (\ref{connection-second order-x}) holds.

\vspace{1mm}

{\bf Remark 3.2}\ Concerning that $(\bar{X}^{t,x;\bar{u}}(\cdot),\bar{u}(\cdot))=(0,\bar{u}(\cdot))$ is an optimal pair when $x=0$, one can also check it by the stochastic verification theorem of \cite{Zhang12} with test function $\varphi(x)=-x$. We mention that the non-smooth version of stochastic verification theorem for the non-recursive optimal control problem was first studied in Zhou, Yong and Li \cite{ZYL97}, and was correctly proved in Gozzi, \'{S}wie\c{c}h and Zhou \cite{GSZ05,GSZ10}.

\vspace{1mm}

\begin{theorem}\label{theorem-NSW-ACC2016}
Let {\bf (H1), (H2), (H3)} hold
and $(t,x)\in[0,T)\times\mathbf{R}^n$ be fixed. Suppose that $\bar{u}(\cdot)$ is an optimal control for {\bf Problem (SROCP)}, and
$(\bar{X}^{t,x;\bar{u}}(\cdot),\bar{Y}^{t,x;\bar{u}}(\cdot),\\\bar{Z}^{t,x;\bar{u}}(\cdot))$
is the optimal trajectory. Let $(p^*(\cdot),q^*(\cdot),k^*(\cdot))\in S^2_\mathcal{F}([0,T];\mathbf{R}^n)\times S^2_\mathcal{F}([0,T];\\\mathbf{R})\times L^2_\mathcal{F}([0,T];\mathbf{R}^n)$ satisfy the following first-order adjoint FBSDE (\cite{Peng93})
\begin{equation}\label{adjoint first-order FBSDE}
\left\{
\begin{aligned}
-dp^*(s)&=\big[\bar{b}_x(s)^\top p^*(s)-\bar{f}_x(s)^\top q^*(s)+\bar{\sigma}_x(s)k^*(s)\big]ds-k^*(s)dW(s),\\
 dq^*(s)&=\bar{f}_y(s) q^*(s)ds+\bar{f}_z(s) q^*(s)dW(s),\ s\in[t,T],\\
  p^*(T)&=-\phi_x(\bar{X}^{t,x;\bar{u}}(T)) q^*(T),\quad q^*(t)=1,
\end{aligned}
\right.
\end{equation}
then
\begin{equation}\label{conection-first order x}
\begin{aligned}
D_x^{1,-}V(s,\bar{X}^{t,x;\bar{u}}(s))\subset\{p^*(s)q^*(s)^{-1}\}\subset D_x^{1,+}V(s,\bar{X}^{t,x;\bar{u}}(s)),
\end{aligned}
\end{equation}
for all $s\in[t,T],\mathbf{P}\mbox{-}a.s.$, where $V(\cdot,\cdot)$ is the value function defined by (\ref{value function}).
\end{theorem}

\vspace{1mm}

{\bf Remark 3.3}\ This result has been proved in \cite{NSW16} with the convex control domain. By comparing Theorem \ref{main theorem-x} with Theorem \ref{theorem-NSW-ACC2016}, it is not hard to guess that $p(s)=-p^*(s)q^*(s)^{-1}$, which is in fact true as shown in the following proof. This means that the first-order adjoint equation (\ref{adjoint equation-first order}) introduced by \cite{Hu15} can be derived directly from the adjoint FBSDE of \cite{Peng93}. Conversely, if we can show that $p(s)=-p^*(s)q^*(s)^{-1}$, then Theorem \ref{theorem-NSW-ACC2016} is indeed a special case of Theorem \ref{main theorem-x} if we only consider the first-order case. Let us give the detailed proof as follows.

\vspace{1mm}

{\it Proof of Theorem \ref{theorem-NSW-ACC2016}}\quad We only need to prove that $p(s)=-p^*(s)q^*(s)^{-1}$, for all $s\in[t,T]$. In fact, from (\ref{adjoint first-order FBSDE}), we obtain
$q^*(s)=\exp\Big\{\int_t^s\bar{f}_y(r)dr-\frac{1}{2}\int_t^s|\bar{f}_z(r)|^2dr+\int_t^s\bar{f}_z(r)dW(r)\Big\}$. Thus applying It\^{o}'s formula to $p^*(s)q^*(s)^{-1}$, we get
\begin{equation*}
\begin{aligned}
&d\{p^*(s)q^*(s)^{-1}\}=\Big\{-\bar{f}_y(s) p^*(s)q^*(s)^{-1}-\bar{b}_x(s)^\top p^*(s)q^*(s)^{-1}\\
&\quad+|\bar{f}_z(s)|^2p^*(s)q^*(s)^{-1}-\bar{\sigma}_x(s)^\top k^*(s)q^*(s)^{-1}-\bar{f}_z(s) k^*(s)q^*(s)^{-1}+\bar{f}_x(s)\Big\}ds\\
&\quad+q^*(s)^{-1}\big[k^*(s)-p^*(s)\bar{f}_z(s)^\top\big]dW(s).
\end{aligned}
\end{equation*}
Let $\mathfrak{p}(s)=-p^*(s)q^*(s)^{-1},\mathfrak{q}(s)=-q^*(s)^{-1}\big[k^*(s)-p^*(s)\bar{f}_z(s)^\top\big]$, we have
\begin{equation*}
\begin{aligned}
-d\mathfrak{p}(s)
      &=\Big\{\bar{f}_y(s)\mathfrak{p}(s)+\bar{f}_z(s)\bar{\sigma}_x(s)^\top\mathfrak{ p}(s)+\bar{b}_x(s)^\top \mathfrak{p}(s)+\bar{f}_z(s)\mathfrak{q}(s)\\
      &\qquad+\bar{\sigma}_x(s)^\top \mathfrak{q}(s)+\bar{f}_x(s)\Big\}ds-\mathfrak{q}(s)dW(s).
\end{aligned}
\end{equation*}
Comparing this with (\ref{adjoint equation-first order}) and noticing that $\mathfrak{p}(T)=-p^*(T)q^*(T)^{-1}=\phi_x(\bar{X}^{t,x;\bar{u}}(T))$, we get $p(s)=\mathfrak{p}(s)=-p^*(s)q^*(s)^{-1}$ from the uniqueness of the solution to BSDE (\ref{adjoint equation-first order}). The proof is complete. \quad$\Box$

\subsection{Differentials in the time variable}

In this subsection, we extend (\ref{Shi10-2}) to the case with non-convex control domain and non-smooth value function.

Let us recall the notion of the right super- and sub-jets in the time variable $t$, see \cite{YZ99}. For $v\in C([0,T]\times\mathbb{R}^n)$, and $(\hat{t},\hat{x})\in[0,T]\times\mathbb{R}^n$, we define
\begin{equation}\label{right super- and sub-jets time}
\left\{
\begin{aligned}
D_{t+}^{1,+}v(\hat{t},\hat{x})&:=\Big\{q\in\mathbf{R}\big|v(t,\hat{x})\leq v(\hat{t},\hat{x})+q(t-\hat{t})
                               +o(|t-\hat{t}|),\mbox{ as }t\downarrow\hat{t}\Big\},\\
D_{t+}^{1,-}v(\hat{t},\hat{x})&:=\Big\{q\in\mathbf{R}\big|v(t,\hat{x})\geq v(\hat{t},\hat{x})+q(t-\hat{t})
                               +o(|t-\hat{t}|),\mbox{ as }t\downarrow\hat{t}\Big\}.
\end{aligned}
\right.
\end{equation}

\vspace{1mm}

\begin{theorem}\label{main theorem-t}
Let {\bf (H1), (H2), (H3)} hold
and $(t,x)\in[0,T)\times\mathbf{R}^n$ be fixed. Suppose that $\bar{u}(\cdot)$ is an optimal control for {\bf Problem (RSOCP)}, and
$(\bar{X}^{t,x;\bar{u}}(\cdot),\bar{Y}^{t,x;\bar{u}}(\cdot),\\\bar{Z}^{t,x;\bar{u}}(\cdot))$
is the optimal state. Let $(p(\cdot),q(\cdot))\in S^2_\mathcal{F}([t,T];\mathbf{R}^n)\times L^2_\mathcal{F}([t,T];\mathbf{R}^n)$ and $(P(\cdot),Q(\cdot))\in S^2_\mathcal{F}([t,T];\mathcal{S}^n)\times L^2_\mathcal{F}([t,T];\mathcal{S}^n)$ satisfy the first- and second-order adjoint equation (\ref{adjoint equation-first order}) and (\ref{adjoint equation-second order}), respectively. Let $V(\cdot,\cdot)$ be the value function. Then
\begin{equation}\label{connection-time}
\begin{aligned}
\big[\mathcal{H}_1\big(s,\bar{X}^{t,x;\bar{u}}(s),\bar{u}(s)\big),\infty\big)\subseteq D_{t+}^{1,+}V(s,\bar{X}^{t,x;\bar{u}}(s)),\ a.e.\ s\in[t,T],\mathbf{P}\mbox{-}a.s.,
\end{aligned}
\end{equation}
\begin{equation}\label{connection-time 2}
\begin{aligned}
D_{t+}^{1,-}V(s,\bar{X}^{t,x;\bar{u}}(s))\subseteq \big(-\infty,\mathcal{H}_1\big(s,\bar{X}^{t,x;\bar{u}}(s),\bar{u}(s)\big)\big],\ a.e.\ s\in[t,T],\mathbf{P}\mbox{-}a.s.,
\end{aligned}
\end{equation}
where $\mathcal{H}_1:[0,T]\times\mathbf{R}^n\times\mathbf{U}\rightarrow\mathbf{R}$ is defined as
\begin{equation}\label{mathcal H1}
{\small\begin{aligned}
\mathcal{H}_1(t,x,u)&:=\mathcal{H}\big(t,x,-V(t,x),\bar{\sigma}(t)^\top p(t),p(t),q(t),P(t),u\big)-\frac{1}{2}\mbox{tr}\big\{\bar{\sigma}(t)^\top P(t)\bar{\sigma}(t)\big\}\\
&=f\big(t,x,-V(t,x),\sigma(t,x,u)^\top p(t),u\big)+\langle p(t),b(t,x,u)\rangle\\
&\quad\ +\big\langle q(t)-P(t)\bar{\sigma}(t,x,u),\sigma(t,x,u)\big\rangle+\frac{1}{2}\big\langle P(t)\sigma(t,x,u),\sigma(t,x,u)\big\rangle\\
&=G\big(t,x,-V(t,x),p(t),P(t),u\big)+\big\langle q(t)-P(t)\sigma(t,\bar{x},\bar{u}),\sigma(t,x,u)\big\rangle.
\end{aligned}}
\end{equation}
\end{theorem}
We postpone the proof of Theorem \ref{main theorem-t} to the next section.

\vspace{1mm}

{\bf Remark 3.4}\ Similar to the classical result (see pp. 267 in \cite{YZ99}), in Theorem \ref{main theorem-t} we only have the {\it right} super- (resp. sub-) jets and it can not be extended to the two-sided super- (sub-) jets. This is due to the adaptiveness requirement.

\vspace{1mm}

{\bf Example 3.1 (continued)}\quad Consider the problem in Example 3.1 again. From (\ref{value function-exapmle 1}), we have for $x=0$,
\begin{equation*}
D_{t+}^{1,+}V(s,\bar{X}^{t,x;\bar{u}}(s))=[0,\infty),\quad D_{t+}^{1,-}V(s,\bar{X}^{t,x;\bar{u}}(s))=(-\infty,0].
\end{equation*}
Since for arbitrary $\bar{u}(\cdot)$ we have $\bar{X}^{t,x;\bar{u}}(s)=0$ and $q(s)=P(s)=Q(s)=0$ for all $s\in[t,T]$, which implies that $\mathcal{H}_1(s,\bar{X}^{t,x;\bar{u}}(s),\bar{u}(s))=0$, $s\in[t,T]$. Thus (\ref{connection-time}) and (\ref{connection-time 2}) hold.
For $x\leq0$, we have
\[D_{t+}^{1,+}V(s,\bar{X}^{t,x;\bar{u}}(s)=[V_s(s,\bar{X}^{t,x;\bar{u}}(s)),\infty)
=[-e^{s-T}\bar{X}^{t,x;\bar{u}}(s),\infty),\]
and
\[
D_{t+}^{1,-}V(s,\bar{X}^{t,x;\bar{u}}(s))=
(-\infty,V_s(s,\bar{X}^{t,x;\bar{u}}(s))]=(-\infty,-e^{s-T}\bar{X}^{t,x;\bar{u}}(s)].\]
Since now $\bar{u}(s)=-1$ or $2$, $p(s)=e^{s-T}$, $q(s)=P(s)=Q(s)=0$ for all $s\in[t,T]$, thus
\[
\begin{aligned}
&\mathcal{H}_1(s,\bar{X}^{t,x;\bar{u}}(s),\bar{u}(s))=p(s)[\bar{X}^{t,x;\bar{u}}(s)+\bar{X}^{t,x;\bar{u}}(s)\bar{u}(s)]-\bar{X}^{t,x;\bar{u}}(s)p(s)(\bar{u}(s))^{2}\\
=&-e^{s-T}\bar{X}^{t,x;\bar{u}}(s).
\end{aligned}
\]
Thus (\ref{connection-time}) and (\ref{connection-time 2}) hold.
For $x\ge0$, we have
\[D_{t+}^{1,+}V(s,\bar{X}^{t,x;\bar{u}}(s))=[V_s(s,\bar{X}^{t,x;\bar{u}}(s)),\infty)
=[e^{T-s}\bar{X}^{t,x;\bar{u}}(s),\infty),\]
and
\[D_{t+}^{1,-}V(s,\bar{X}^{t,x;\bar{u}}(s))=(-\infty,V_s(s,\bar{X}^{t,x;\bar{u}}(s))]
=(-\infty,e^{T-s}\bar{X}^{t,x;\bar{u}}(s)].\]
Since now $\bar{u}(s)=0$ or $1$, $p(s)=e^{T-s}$, $q(s)=P(s)=Q(s)=0$ for $s\in[t,T]$, thus
\[
\begin{aligned}
&\mathcal{H}_1(s,\bar{X}^{t,x;\bar{u}}(s),\bar{u}(s))=p(s)[\bar{X}^{t,x;\bar{u}}(s)+\bar{X}^{t,x;\bar{u}}(s)\bar{u}(s)]-\bar{X}^{t,x;\bar{u}}(s)p(s)(\bar{u}(s))^{2}\\
=&\ e^{T-s}\bar{X}^{t,x;\bar{u}}(s).
\end{aligned}
\]
Thus (\ref{connection-time}) and (\ref{connection-time 2}) hold.

%
\subsection{Smooth Case}

In this subsection, we will give the smooth version of Theorem \ref{main theorem-x} and Theorem \ref{main theorem-t}. Different from \cite{Shi10,SY13}, here the control domain can be non-convex. It is also worth to emphasize that we can derive the first- and second-order adjoint equations (\ref{adjoint equation-first order}), (\ref{adjoint equation-second order}) of \cite{Hu15} directly, see Remark 3.6.

\vspace{1mm}

\begin{corollary}\label{corollary-smooth case}
Let {\bf (H1), (H2), (H3)} hold
and $(t,x)\in[0,T)\times\mathbf{R}^n$ be fixed. Suppose that $\bar{u}(\cdot)$ is an optimal control for {\bf Problem (SROCP)}, and
$(\bar{X}^{t,x;\bar{u}}(\cdot),\\\bar{Y}^{t,x;\bar{u}}(\cdot),\bar{Z}^{t,x;\bar{u}}(\cdot))$
is the optimal trajectory. If $V(\cdot,\cdot)\in C^{1,2}([0,T]\times\mathbf{R}^n)$, then
\begin{equation}\label{relation with respect to time variable}
\begin{aligned}
  &V_s(s,\bar{X}^{t,x;\bar{u}}(s))\\
 =&\ G\big(s,\bar{X}^{t,x;\bar{u}}(s),-V(s,\bar{X}^{t,x;\bar{u}}(s)),-V_x(s,\bar{X}^{t,x;\bar{u}}(s)),-V_{xx}(s,\bar{X}^{t,x;\bar{u}}(s)),\bar{u}(s)\big)\\
 =&\ \max\limits_{u\in\mathbf{U}}G\big(s,\bar{X}^{t,x;\bar{u}}(s),-V(s,\bar{X}^{t,x;\bar{u}}(s)),
 -V_x(s,\bar{X}^{t,x;\bar{u}}(s)),-V_{xx}(s,\bar{X}^{t,x;\bar{u}}(s)),u\big),
\end{aligned}
\end{equation}
a.e.$s\in[t,T],\mathbf{P}\mbox{-}$a.s., where $G$ is defined as (\ref{generalized Hamiltonian}). Moreover, if $V(\cdot,\cdot)\in C^{1,3}([0,T]\times\mathbf{R}^n)$ and $V_{tx}(\cdot,\cdot)$ is continuous, then
\begin{equation}\label{relation Hu and Peng-p,q}
\left\{
\begin{aligned}
  p(s)=&-V_x(s,\bar{X}^{t,x;\bar{u}}(s)),\qquad\qquad\text{ for all } s\in[t,T],\ \mathbf{P}\mbox{-}a.s.,\\
  q(s)=&-V_{xx}(s,\bar{X}^{t,x;\bar{u}}(s))\bar{\sigma}(s),\qquad\quad a.e.\ s\in[t,T],\ \mathbf{P}\mbox{-}a.s.,\\
\end{aligned}
\right.
\end{equation}
where $(p(\cdot),q(\cdot))$ satisfy (\ref{adjoint equation-first order}). Furthermore, if $V(\cdot,\cdot)\in C^{1,4}([0,T]\times\mathbf{R}^n)$ and $V_{txx}(\cdot,\cdot)$ is continuous, then
\begin{equation}\label{relation Hu and Peng-P,Q}
\begin{aligned}
  -V_{xx}(s,\bar{X}^{t,x;\bar{u}}(s))\geq P(s),\qquad\text{ for all } s\in[t,T],\ \mathbf{P}\mbox{-}a.s.,
\end{aligned}
\end{equation}
where $(P(\cdot),Q(\cdot))$ satisfies (\ref{adjoint equation-second order}).
\end{corollary}

\vspace{1mm}

{\bf Remark 3.5}\ Comparing (\ref{relation Hu and Peng-p,q}) with (\ref{Shi10}), similar to the non-smooth case, we can obtain also that \begin{equation}\label{relation Hu and Peng-smooth case}
\left\{
\begin{aligned}
p(s)&=-p^*(s)q^*(s)^{-1},\\
q(s)&=-q^*(s)^{-1}\big[k^*(s)-p^*(s)\bar{f}_z(s)^\top\big].
\end{aligned}
\right.
\end{equation}
This relation also reveals the connection between $(p,q)$ (solution to BSDE (\ref{adjoint equation-first order})) of \cite{Hu15} with $(p^*,q^*,k^*)$ (solution to FBSDE (\ref{adjoint first-order FBSDE})) of \cite{Peng93}.

\vspace{1mm}

{\it Proof of Corollary \ref{corollary-smooth case}.}\quad For fixed $t\in[0,T)$, by the backward semigroup property of \cite{Peng97}, since $\bar{u}(\cdot)$ is the optimal control, we obtain $-V(s,\bar{X}^{t,x;\bar{u}}(s))=\bar{Y}^{t,x;\bar{u}}(s), s\in[t,T]$ which satisfies BSDE (\ref{controlled BSDE}). Applying It\^{o}'s formula to
$V(s,\bar{X}^{t,x;\bar{u}}(s))$, we have
\begin{equation*}\label{dV(t)}
\begin{aligned}
 &dV(s,\bar{X}^{t,x;\bar{u}}(s))=\big[V_s(s,\bar{X}^{t,x;\bar{u}}(s))+\big\langle V_x(s,\bar{X}^{t,x;\bar{u}}(s)),\bar{b}(s)\big\rangle\\
&\ +\frac{1}{2}\mbox{tr}\big(\bar{\sigma}(s)^\top V_{xx}(s,\bar{X}^{t,x;\bar{u}}(s))\bar{\sigma}(s)\big)\big]ds+V_x(s,\bar{X}^{t,x;\bar{u}}(s))^\top\bar{\sigma}(s)dW(s).
\end{aligned}
\end{equation*}
Comparing this with BSDE (\ref{controlled BSDE}), we conclude that
\begin{equation}\label{1 equation}
\left\{
\begin{aligned}
&V_s(s,\bar{X}^{t,x;\bar{u}}(s))+\big\langle V_x(s,\bar{X}^{t,x;\bar{u}}(s)),\bar{b}(s)\big\rangle
+\frac{1}{2}\mbox{tr}\big(\bar{\sigma}(s)^\top V_{xx}(s,\bar{X}^{t,x;\bar{u}}(s))\bar{\sigma}(s)\big)\\
&=\bar{f}(s)=f\big(s,\bar{X}^{t,x;\bar{u}}(s),\bar{Y}^{t,x;\bar{u}}(s),\bar{Z}^{t,x;\bar{u}}(s),\bar{u}(s)\big),\text{ for all } s\in[t,T],\ \mathbf{P}\mbox{-}a.s.,\\
&V_x(s,\bar{X}^{t,x;\bar{u}}(s))^\top\bar{\sigma}(s)=-\bar{Z}^{t,x;\bar{u}}(s),\qquad\quad\qquad a.e.s\in[t,T],\mathbf{P}\mbox{-}a.s.
\end{aligned}
\right.
\end{equation}
It follows from (\ref{1 equation}) that the first equality of (\ref{relation with respect to time variable}) holds. Moreover,
since $V(\cdot,\cdot)\in C^{1,2}([0,T]\times\mathbf{R}^n)$, it is the classical solution to (\ref{HJB equation}), thus using Proposition \ref{proposition-Peng1997} we obtain the second equality of (\ref{relation with respect to time variable}).

On the other hand, by (\ref{HJB equation}) and  (\ref{relation with respect to time variable}), it follows that for all $x\in\mathbf{R}^n$,
\begin{equation}\label{HJB equation-result}
\begin{aligned}
   0=&\ -V_s(s,\bar{X}^{t,x;\bar{u}}(s))\\
     &\ +G\big(s,\bar{X}^{t,x;\bar{u}}(s),-V(s,\bar{X}^{t,x;\bar{u}}(s)),-V_x(s,\bar{X}^{t,x;\bar{u}}(s)),-V_{xx}(s,\bar{X}^{t,x;\bar{u}}(s)),\bar{u}(s)\big)\\
 \geq&-V_s(s,x)+G\big(s,x,-V(s,x),-V_x(s,x),-V_{xx}(s,x),\bar{u}(s)\big).
\end{aligned}
\end{equation}
Consequently, if $V(\cdot,\cdot)\in C^{1,3}([0,T]\times\mathbf{R}^n)$, then for all $s\in[t,T]$,
\begin{equation}\label{first-order maximum conditon}
\frac{\partial}{\partial x}\Big\{-V_s(s,x)+G\big(s,x,-V(s,x),-V_x(s,x),-V_{xx}(s,x),\bar{u}(s)\big)\Big\}\Big|_{x=\bar{X}^{t,x;\bar{u}}(s)}=0,
\end{equation}
which is the first-order maximum condition. Furthermore, if $V(\cdot,\cdot)\in C^{1,4}([0,T]\times\mathbf{R}^n)$, the following second-order maximum condition holds: For all $s\in[t,T]$,
\begin{equation}\label{second-order maximum conditon}
\frac{\partial^2}{\partial x^2}\Big\{-V_s(s,x)+G\big(s,x,-V(s,x),-V_x(s,x),-V_{xx}(s,x),\bar{u}(s)\big)\Big\}\Big|_{x=\bar{X}^{t,x;\bar{u}}(s)}\leq0,\
\end{equation}
On the one hand, (\ref{first-order maximum conditon}) yields that (recall (\ref{generalized Hamiltonian})) for all $s\in[t,T]$,
\begin{equation}\label{2 equation}
\begin{aligned}
0&=-V_{sx}(s,\bar{X}^{t,x;\bar{u}}(s))-V_{xx}(s,\bar{X}^{t,x;\bar{u}}(s))\bar{b}(s)-\bar{b}_x(s)^\top V_x(s,\bar{X}^{t,x;\bar{u}}(s))\\
 &\quad-\frac{1}{2}\mbox{tr}\big(\bar{\sigma}(s)^\top V_{xxx}(s,\bar{X}^{t,x;\bar{u}}(s))\bar{\sigma}(s)\big)-\bar{\sigma}_x(s)^\top V_{xx}(s,\bar{X}^{t,x;\bar{u}}(s))\bar{\sigma}(s)+\bar{f}_x(s)\\
 &\quad-\bar{f}_y(s)V_x(s,\bar{X}^{t,x;\bar{u}}(s))
 -\bar{f}_z(s)\big[V_{xx}(s,\bar{X}^{t,x;\bar{u}}(s))\bar{\sigma}(s)+\bar{\sigma}_x(s)^\top V_x(s,\bar{X}^{t,x;\bar{u}}(s))\big],
\end{aligned}
\end{equation}
where
$\mbox{tr}\big(\bar{\sigma}^\top V_{xxx}\bar{\sigma}\big):=\big(\mbox{tr}\big(\bar{\sigma}^\top\big((V_x)^1\big)_{xx}\bar{\sigma}\big),\cdots,
\mbox{tr}\big(\bar{\sigma}^\top\big((V_x)^n\big)_{xx}\bar{\sigma}\big)\big)^\top$, with $\big((V_x)^1,\\\cdots,(V_x)^n\big)^\top=V_x$. Applying It\^{o}'s formula to $V_x(s,\bar{X}^{t,x;\bar{u}}(s))$, using (\ref{2 equation}) and the fact that $V_{sx}(\cdot,\cdot)$ is continuous we get
\begin{equation*}
\begin{aligned}
 &-dV_x(s,\bar{X}^{t,x;\bar{u}}(s))=\Big\{\bar{b}_x(s)^\top V_x(s,\bar{X}^{t,x;\bar{u}}(s))+\bar{\sigma}_x(s)^\top V_{xx}(s,\bar{X}^{t,x;\bar{u}}(s))\bar{\sigma}(s)-\bar{f}_x(s)\\
 &\ +\bar{f}_y(s)V_x(s,\bar{X}^{t,x;\bar{u}}(s))+\bar{f}_z(s)\big[V_{xx}(s,\bar{X}^{t,x;\bar{u}}(s))\bar{\sigma}(s)+\bar{\sigma}_x(s)^\top V_x(s,\bar{X}^{t,x;\bar{u}}(s))\big]\Big\}ds\\
 &\ -V_{xx}(s,\bar{X}^{t,x;\bar{u}}(s))\bar{\sigma}(s)dW(s).
\end{aligned}
\end{equation*}
Note that $V$ solves (\ref{HJB equation}), thus $V_x(T,\bar{X}^{t,x;\bar{u}}(T))=-\phi_x(\bar{X}^{t,x;\bar{u}}(T))$. Then by the uniqueness of the solutions to (\ref{adjoint equation-first order}), we obtain (\ref{relation Hu and Peng-p,q}).

Moreover, (\ref{second-order maximum conditon}) yields that, for any $s\in[t,T]$,
\begin{equation}\label{comparison}
{\small\begin{aligned}
 &V_{sxx}(s,\bar{X}^{t,x;\bar{u}}(s))+V_{xxx}(s,\bar{X}^{t,x;\bar{u}}(s))\bar{b}(s)
  +\frac{1}{2}\big[\bar{\sigma}(s)^\top V_{xxxx}(s,\bar{X}^{t,x;\bar{u}}(s))\bar{\sigma}(s)\big]\\
 &\geq-V_{xx}(s,\bar{X}^{t,x;\bar{u}}(s))\bar{b}_x(s)-\bar{b}_x(s)^\top V_{xx}(s,\bar{X}^{t,x;\bar{u}}(s))-\bar{f}_y(s)V_{xx}(s,\bar{X}^{t,x;\bar{u}}(s))\\
 &\quad-\bar{f}_z(s)\bar{\sigma}_x(s)^\top V_{xx}(s,\bar{X}^{t,x;\bar{u}}(s))-\bar{f}_z(s)V_{xx}(s,\bar{X}^{t,x;\bar{u}}(s))\bar{\sigma}_x(s)\\
 &\quad-\bar{\sigma}_x(s)^\top V_{xx}(s,\bar{X}^{t,x;\bar{u}}(s))\bar{\sigma}_x(s)-V_{xxx}(s,\bar{X}^{t,x;\bar{u}}(s))\bar{\sigma}(s)\bar{\sigma}_x(s)\\
 &\quad-\bar{f}_z(s)V_{xxx}(s,\bar{X}^{t,x;\bar{u}}(s))\bar{\sigma}(s)-\bar{\sigma}_x(s)^\top V_{xxx}(s,\bar{X}^{t,x;\bar{u}}(s))\bar{\sigma}(s)\\
 &\quad-\bar{b}_{xx}(s)^\top V_x(s,\bar{X}^{t,x;\bar{u}}(s))-\bar{\sigma}_{xx}(s)^\top \bar{f}_z(s)V_x(s,\bar{X}^{t,x;\bar{u}}(s))-\bar{\sigma}_{xx}(s)^\top V_{xx}(s,\bar{X}^{t,x;\bar{u}}(s))\bar{\sigma}(s)\\
 &\quad+\big[I_{n\times n},-V_x(s,\bar{X}^{t,x;\bar{u}}(s)),-\bar{\sigma}_x(s)^\top V_x(s,\bar{X}^{t,x;\bar{u}}(s))-V_{xx}(s,\bar{X}^{t,x;\bar{u}}(s))\bar{\sigma}(s)\big]D^2\bar{f}(s)\\
 &\qquad\cdot\big[I_{n\times n},-V_x(s,\bar{X}^{t,x;\bar{u}}(s)),-\bar{\sigma}_x(s)^\top V_x(s,\bar{X}^{t,x;\bar{u}}(s))-V_{xx}(s,\bar{X}^{t,x;\bar{u}}(s))\bar{\sigma}(s)\big]^\top.
\end{aligned}}
\end{equation}
In the above and sequel, the notations of partial derivatives have their own definitions which we will not clarify on by one, for the limit of the space (for simplicity, the readers can verify the calculus just using $n=1$, i.e. $x$ is one dimensional).

Applying It\^{o}'s formula to $V_{xx}(s,\bar{X}^{t,x;\bar{u}}(s))$, we obtain
\begin{equation*}
\begin{aligned}
 &dV_{xx}(s,\bar{X}^{t,x;\bar{u}}(s))=\Big\{V_{xxs}(s,\bar{X}^{t,x;\bar{u}}(s))+V_{xxx}(s,\bar{X}^{t,x;\bar{u}}(s))\bar{b}(s)\\
 &\quad+\frac{1}{2}\big[\bar{\sigma}(s)^\top V_{xxxx}(s,\bar{X}^{t,x;\bar{u}}(s))\bar{\sigma}(s)\big]\Big\}ds-V_{xxx}(s,\bar{X}^{t,x;\bar{u}}(s))^\top\bar{\sigma}(s)dW(s).
\end{aligned}
\end{equation*}
Define $\mathcal{P}(s)=-V_{xx}(s,\bar{X}^{t,x;\bar{u}}(s)),\mathcal{Q}(s)=-V_{xxx}(s,\bar{X}^{t,x;\bar{u}}(s))\bar{\sigma}(s)$, we have
\begin{equation}\label{comparied matrix valued BSDE}
\begin{aligned}
 -d\mathcal{P}(s)=&\Big\{V_{xxs}(s,\bar{X}^{t,x;\bar{u}}(s))+V_{xxx}(s,\bar{X}^{t,x;\bar{u}}(s))\bar{b}(s)\\
                  &+\frac{1}{2}\big[\bar{\sigma}(s)^\top V_{xxxx}(s,\bar{X}^{t,x;\bar{u}}(s))\bar{\sigma}(s)\big]\Big\}ds-\mathcal{Q}(s)dW(s).
\end{aligned}
\end{equation}
From (\ref{comparison}) and the continuity of $V_{sxx}(\cdot,\cdot)$ as well as (\ref{relation Hu and Peng-p,q}), we have
\begin{equation*}\label{comparison-1}
{\small\begin{aligned}
 &V_{sxx}(s,\bar{X}^{t,x;\bar{u}}(s))+V_{xxx}(s,\bar{X}^{t,x;\bar{u}}(s))\bar{b}(s)
  +\frac{1}{2}\big[\bar{\sigma}(s)^\top V_{xxxx}(s,\bar{X}^{t,x;\bar{u}}(s))\bar{\sigma}(s)\big]\\
 &\geq\bar{f}_y(s)\mathcal{P}(s)+\big[\bar{f}_z(s)\bar{\sigma}_x(s)^\top+\bar{b}_x(s)^\top\big]\mathcal{P}(s)
  +\mathcal{P}(s)\big[\bar{f}_z(s)\bar{\sigma}_x(s)+\bar{b}_x(s)\big]+\bar{\sigma}_x(s)^\top\mathcal{P}(s)\bar{\sigma}_x(s)\\
 &\quad+\bar{f}_z(s)\mathcal{Q}(s)+\bar{\sigma}_x(s)^\top\mathcal{Q}(s)+\mathcal{Q}(s)\bar{\sigma}_x(s)+\bar{b}_{xx}(s)^\top p(s)+\bar{\sigma}_{xx}(s)^\top\big[\bar{f}_z(s)p(s)\\
 &\quad+q(s)\big]+\big[I_{n\times n},p(s),\bar{\sigma}_x(s)^\top p(s)+q(s)\big]D^2\bar{f}(s)\big[I_{n\times n},p(s),\bar{\sigma}_x(s)^\top p(s)+q(s)\big]^\top.
\end{aligned}}
\end{equation*}
Note that $P(T)=\mathcal{P}(T)=-V_{xx}(T,\bar{X}^{t,x;\bar{u}}(T))=\phi_{xx}(\bar{X}^{t,x;\bar{u}}(T))$. Using above relation, one can check the condition (condition (13) of Hu and Peng \cite{HuPeng06}) of the comparison theorem for the matrix-valued BSDEs  (\ref{adjoint equation-second order}) and (\ref{comparied matrix valued BSDE}) holds, thus we obtain (\ref{relation Hu and Peng-P,Q}). The proof is complete.\quad$\Box$

\vspace{1mm}

{\bf Remark 3.6}\ From the above proof, we can observe that if the second-order derivative of (\ref{second-order maximum conditon}) at $x=\bar{X}^{t,x;\bar{u}}(s)$ equals to zero, then we have $-V_{xx}(s,\bar{X}^{t,x;\bar{u}}(s))=P(s),\text{ for all } s\in[t,T],\mathbf{P}\mbox{-}a.s.$ The following example is in this case.

\vspace{1mm}

{\bf Example 3.2}\quad Consider the following controlled SDE ($n=1$):
\begin{equation}\label{controlled SDE:example 2}
\left\{
\begin{aligned}
 dX^{t,x;u}(s)&=2u(s)ds+u(s)dW(s),\ s\in[t,T],\\
  X^{t,x;u}(t)&=x,
\end{aligned}
\right.
\end{equation}
with $\mathbf{U}=[-3,-2]\bigcup[1,2]$. The cost functional is defined as (\ref{cost functional}) with
\begin{equation}\label{controlled BSDE:example 2}
\left\{
\begin{aligned}
-dY^{t,x;u}(s)&=\Big\{\frac{1}{2}|u(s)\mbox{th}X^{t,x;u}(s)|^2-u(s)\mbox{th}X^{t,x;u}(s)-|u(s)|^2-u(s)\\
              &\qquad-Z^{t,x;u}(s)\Big\}ds-Z^{t,x;u}(s)dW(s),\ s\in[t,T],\\
  Y^{t,x;u}(T)&=\ln \mbox{ch}X^{t,x;u}(T),
\end{aligned}
\right.
\end{equation}
where $\mbox{ch}x:=\frac{1}{2}(e^x+e^{-x})$ and $\mbox{th}x=\frac{e^x-e^{-x}}{e^x+e^{-x}}$.
Then the HJB equation (\ref{HJB equation}) writes
\[
\left\{
\begin{aligned}
&-v_t(t,x)+\sup_{u\in\mathbf{U}}\left\{-\frac{1}{2}v_{xx}(t,x)u^2+\frac{1}{2}|u\mbox{th}x|^2-u\mbox{th}x-v_{x}(t,x)u-u^{2}-u\right\}=0,\\ &\qquad\qquad\qquad\qquad\qquad\qquad\qquad\qquad\qquad\qquad\qquad\qquad\qquad (t,x)\in[0,T]\times\mathbf{R}^n,\\
&v(T,x)=-\ln \mbox{ch}x,\quad \text{ for all } x\in\mathbf{R}^n,
\end{aligned}
\right.
\]
which admits a unique solution $V(t,x)=-\ln \mbox{ch}x$. One can check that $V_t(t,x)\equiv0$, $V_x(t,x)=-\mbox{th}x$ and $V_{xx}(t,x)=-(\mbox{ch}x)^{-2}$.

Let $\bar{u}(\cdot)$ be an optimal control, then the first-order adjoint equation writes
\begin{equation}\label{first order equation:example 2}
\left\{
\begin{aligned}
-dp(s)&=\Big\{|\bar{u}(s)|^2\mbox{th}\bar{X}^{t,x;\bar{u}}(s)(\mbox{ch}\bar{X}^{t,x;\bar{u}}(s))^{-2}-\bar{u}(s)(\mbox{ch}\bar{X}^{t,x;\bar{u}}(s))^{-2}\\
      &\qquad-q(s)\Big\}ds-q(s)dW(s),\quad s\in[t,T],\\
  p(T)&=\mbox{th}\bar{X}^{t,x;\bar{u}}(T)),
\end{aligned}
\right.
\end{equation}
which, by Ito's formula, admits a unique solution $(p(s),q(s))=\big(\mbox{th}\bar{X}^{t,x;\bar{u}}(s),\bar{u}(s)\\(\mbox{ch}\bar{X}^{t,x;\bar{u}}(s))^{-2}\big)$. Therefore (\ref{relation Hu and Peng-p,q}) holds. Now let us focus on the second order adjoint equation
\begin{equation}\label{second order equation:example 2}
\left\{
\begin{aligned}
-dP(s)&=g\big(s,\bar{X}^{t,x;\bar{u}}(s),Q(s),\bar{u}(s)\big)ds-Q(s)dW(s),\quad s\in[t,T],\\
  P(T)&=(\mbox{ch}\bar{X}^{t,x;\bar{u}}(T))^{-2},
\end{aligned}
\right.
\end{equation}
where $g(s,x,Q(s),u):=(\mbox{ch}x)^{-4}\big[2u^2-\frac{u^2-u}{2}e^{2x}-\frac{u^2+u}{2}e^{-2x}\big]-Q(s)$, which, by again Ito's formula, admits a unique solution
$(P(s),Q(s))=\big((\mbox{ch}\bar{X}^{t,x;\bar{u}}(s))^{-2},-2\bar{u}(s)\\(\mbox{ch}\bar{X}^{t,x;\bar{u}}(s))^{-2}\mbox{th}\bar{X}^{t,x;\bar{u}}(s)\big)$.

On the other hand, from (\ref{second-order maximum conditon}), if we define $F(s,y):=-V_s(s,y)+G\big(s,y,-V(s,y),\\-V_y(s,y),-V_{yy}(s,y),\bar{u}(s)\big)$ which reduces to $F(s,y)=-\frac{1}{2}|\bar{u}(s)|^2-\bar{u}(s)$, and thus $F_y(s,y)=0,F_{yy}(s,y)=0$, for any $(s,y)\in[0,T]\times\mathbf{R}^n$. Therefore (\ref{relation Hu and Peng-P,Q}) should hold with the equality. In fact, recalling that $V_{xx}(t,x)=-(\mbox{ch}x)^{-2}$ and $P(s)=(\mbox{ch}\bar{X}^{t,x;\bar{u}}(s))^{-2}$, (\ref{relation Hu and Peng-P,Q}) holds indeed with the equality.

We observe that the above arguments do not use the exact form of the optimal control. By applying the stochastic verification theorem in \cite{SY13}, we obtain that the optimal control is $\bar{u}(s)=-2$.

\vspace{1mm}

{\bf Remark 3.7}\ We point out that the strict inequality $-V_{xx}(s,\bar{X}^{t,x;\bar{u}}(s))>P(s)$ happens if the inequality in (\ref{second-order maximum conditon}) holds strictly. Specially, as in the following example when $n=1$, this is due to the strict comparison theorem of scalar BSDEs, see \cite{Peng97}.

\vspace{1mm}

{\bf Example 3.3}\quad Consider the following controlled SDE ($n=1$):
\begin{equation}\label{controlled SDE example 3}
\left\{
\begin{aligned}
 dX^{t,x;u}(s)&=2u(s)ds+u(s)dW(s),\ s\in[t,T],\\
  X^{t,x;u}(t)&=x,
\end{aligned}
\right.
\end{equation}
with $\mathbf{U}=[-1,1]\bigcup[2,4]$. The cost functional is defined as (\ref{cost functional}) with
\begin{equation}\label{controlled BSDE example 3}
\left\{
\begin{aligned}
-dY^{t,x;u}(s)&=\Big\{\frac{1}{2}|u(s)\mbox{th}X^{t,x;u}(s)|^2-\frac{1}{2}|\mbox{th}X^{t,x;u}(s)|^2-|u(s)|^2\\
              &\qquad-Z^{t,x;u}(s)\Big\}ds-Z^{t,x;u}(s)dW(s),\ s\in[t,T],\\
  Y^{t,x;u}(T)&=\ln \mbox{ch}X^{t,x;u}(T).
\end{aligned}
\right.
\end{equation}
Then the HJB equation (\ref{HJB equation}) writes
\begin{equation*}
\left\{
\begin{aligned}
&\sup_{u\in\mathbf{U}}\Big\{-\frac{1}{2}v_{xx}(t,x)u^2-v_{x}(t,x)u+\frac{1}{2}|u\mbox{th}x|^2-u^{2}\Big\}\\
&\hspace{8mm}-v_t(t,x)-\frac{1}{2}|\mbox{th}x|^2=0,\ (t,x)\in[0,T]\times\mathbf{R}^n,\\
&v(T,x)=-\ln \mbox{ch}x,\quad \text{ for all } x\in\mathbf{R}^n,
\end{aligned}
\right.
\end{equation*}
which admits a unique solution $V(t,x)=-\ln \mbox{ch}x$. Moreover, by applying the stochastic verification theorem in \cite{SY13}, we obtain that $\bar{u}(s)=\mbox{th}\bar{X}^{t,x,\bar{u}}(s)$.

The adjoint equations are
\begin{equation}\label{first order equation example 3}
\left\{
\begin{aligned}
-dp(s)&=\Big\{\mbox{th}\bar{X}^{t,x;\bar{u}}(s)(\mbox{ch}\bar{X}^{t,x;\bar{u}}(s))^{-2}\big[(\bar{u}(s))^2-1\big]-q(s)\Big\}ds\\
      &\quad-q(s)dW(s),\quad s\in[t,T],\\
  p(T)&=\mbox{th}\bar{X}^{t,x;\bar{u}}(T)),
\end{aligned}
\right.
\end{equation}
and
\begin{equation}\label{second order equation example 3}
\left\{
\begin{aligned}
-dP(s)&=g_1\big(s,\bar{X}^{t,x;\bar{u}}(s),Q(s),\bar{u}(s)\big)ds-Q(s)dW(s),\quad s\in[t,T],\\
  P(T)&=(\mbox{ch}\bar{X}^{t,x;\bar{u}}(T))^{-2},
\end{aligned}
\right.
\end{equation}
where $g_1(s,x,Q(s),u):=\big[(\mbox{ch}x)^{-4}-2(\mbox{th}x)^2(\mbox{ch}x)^{-2}\big](u^2-1)-Q(s)$. Applying Ito's formula and using $\bar{u}(s)=\mbox{th}\bar{X}^{t,x,\bar{u}}(s)$, we obtain that
$(p(s),q(s))=\big(\mbox{th}\bar{X}^{t,x;\bar{u}}(s),\\(\mbox{ch}\bar{X}^{t,x;\bar{u}}(s))^{-2}\mbox{th}\bar{X}^{t,x;\bar{u}}(s)\big)$ is the unique solution to (\ref{first order equation example 3}).

Now let us focus on (\ref{second order equation example 3}). From (\ref{second-order maximum conditon}), if we define
$\tilde{F}(s,y):=-V_s(s,y)+G\big(s,y,-V(s,y),-V_y(s,y)-V_{yy}(s,y),\bar{u}(s)\big)$ which reduces to $\tilde{F}(s,y)=-\frac{1}{2}(\mbox{th}y)^2-\frac{1}{2}(\mbox{th}\bar{X}^{t,x,\bar{u}}(s))^2+\mbox{th}y\cdot\mbox{th}\bar{X}^{t,x,\bar{u}}(s)$. We can verify that
$\tilde{F}_y(s,\bar{X}^{t,x;\bar{u}}(s))=0$, and $\tilde{F}_{yy}(s,\bar{X}^{t,x;\bar{u}}(s))=-(\mbox{ch}\bar{X}^{t,x;\bar{u}}(s))^{-4}<0$. Next, since $-V_{xx}(t,x)=(\mbox{ch}x)^{-2}$, we denote
$\mathcal{P}(s):=(\mbox{ch}\bar{X}^{t,x;\bar{u}}(s))^{-2},\ \mathcal{Q}(s):=-2\bar{u}(s)(\mbox{ch}\bar{X}^{t,x;\bar{u}}(s))^{-2}\mbox{th}\bar{X}^{t,x;\bar{u}}(s)$.
Applying Ito's formula to $\mathcal{P}(\cdot)$, we have
\begin{equation*}
\left\{
\begin{aligned}
-d\mathcal{P}(s)&=g_2\big(s,\bar{X}^{t,x;\bar{u}}(s),\bar{u}(s)\big)ds-\mathcal{Q}(s)dW(s),\quad s\in[t,T],\\
  \mathcal{P}(T)&=(\mbox{ch}\bar{X}^{t,x;\bar{u}}(T))^{-2},
\end{aligned}
\right.
\end{equation*}
where $g_2(s,x,u):=(\mbox{ch}x)^{-4}\big[2u^2+\frac{-u^2+2u}{2}e^{2x}+\frac{-u^2-2u}{2}e^{-2x}\big]$. We can check directly that for $\bar{u}(s)=\mbox{th}\bar{X}^{t,x,\bar{u}}(s)$,
$g_2(s,\bar{X}^{t,x;\bar{u}}(s),\bar{u}(s))-g_1(s,\bar{X}^{t,x;\bar{u}}(s),\mathcal{Q}(s),\bar{u}(s))=(\mbox{ch}\bar{X}^{t,x;\bar{u}}(s))^{-4}>0$, and the strictly comparison theorem for scalar BSDEs yields that $-V_{xx}(s,\bar{X}^{t,x;\bar{u}}(s))=(\mbox{ch}\bar{X}^{t,x;\bar{u}}(s))^{-2}=\mathcal{P}(s)>P(s)$.

The above arguments are based on the fact that we know the optimal pair $(\bar{X}^{t,x,\bar{u}}(s),\bar{u}(s))$ satisfies $\bar{u}(s)=\mbox{th}\bar{X}^{t,x,\bar{u}}(s)$. Now we will show that this relation
can be found directly from (\ref{relation Hu and Peng-p,q}) as a candidate optimal control. In fact, the optimal control $\bar{u}(\cdot)$ should keep that (\ref{relation Hu and Peng-p,q}) holds where $(p(\cdot),q(\cdot))$ and $(P(\cdot),Q(\cdot))$ solves (\ref{first order equation example 3}) and (\ref{second order equation example 3}) respectively. Thus we should have $p(s)=\mbox{th}\bar{X}^{t,x;\bar{u}}(s)$ and $P(s)\leq(\mbox{ch}\bar{X}^{t,x;\bar{u}}(s))^{-2}$. Applying Ito's formula to $\tilde{p}(s):=\mbox{th}\bar{X}^{t,x;\bar{u}}(s)$, we have
$$
-d\tilde{p}(s)=\tilde{g}(s,\bar{X}^{t,x;\bar{u}}(s),\bar{u}(s))ds-\tilde{q}(s)dW(s),
$$
where $\tilde{q}(s):=\bar{u}(s)(\mbox{ch}\bar{X}^{t,x;\bar{u}}(s))^{-2}$ and $\tilde{g}(s,x,u):=u^2(\mbox{ch}x)^{-2}\mbox{th}x-2u(\mbox{ch}x)^{-2}$. This should coincides with BSDE (\ref{first order equation example 3}). Thus $q(s)=\tilde{q}(s)$ and
$$
[(\bar{u}(s))^2-1]\cdot\big[\mbox{th}x\cdot(\mbox{ch}x)^{-2}\big]-q(s)=\big[(\bar{u}(s))^2(\mbox{ch}x)^{-2}\mbox{th}x-2\bar{u}(s)(\mbox{ch}x)^{-2}\big]\big|_{x=\bar{X}^{t,x;\bar{u}}(s)},
$$
which yields that $\bar{u}(s)=\mbox{th}\bar{X}^{t,x,\bar{u}}(s)$. Moreover, we know that for $\bar{u}(s)=\mbox{th}\bar{X}^{t,x,\bar{u}}(s)$, $\mathcal{P}(s)=(\mbox{ch}\bar{X}^{t,x;\bar{u}}(s))^{-2}>P(s)$ holds. Thus $\bar{u}(s)=\mbox{th}\bar{X}^{t,x,\bar{u}}(s)$ keeps that (\ref{relation Hu and Peng-p,q}) holds. Then it is a candidate optimal control.

Under the assumption that the value function is smooth enough, the following result tells us that the MP can be derived from the DPP, which is an important issue in stochastic control theory.

\vspace{1mm}

\begin{corollary}\label{corollary-DPP to MP}
Let {\bf (H1), (H2), (H3)} hold
and $(t,x)\in[0,T)\times\mathbf{R}^n$ be fixed. Suppose that $\bar{u}(\cdot)$ is an admissible control for {\bf Problem (SROCP)}, and
$(\bar{X}^{t,x;\bar{u}}(\cdot),\bar{Y}^{t,x;\bar{u}}(\cdot),\bar{Z}^{t,x;\bar{u}}(\cdot))$ is the corresponding state. Let $G,\mathcal{H}$ be defined by (\ref{generalized Hamiltonian}), (\ref{Hamiltonian}), respectively, and $(p(\cdot),q(\cdot)),(P(\cdot),Q(\cdot))$ satisfy (\ref{adjoint equation-first order}) and (\ref{adjoint equation-second order}) respectively. If the value function $V(\cdot,\cdot)$ is smooth enough and its derivatives are continuous, and
\begin{equation}\label{DPP}
\begin{aligned}
  &G\big(s,\bar{X}^{t,x;\bar{u}}(s),-V(s,\bar{X}^{t,x;\bar{u}}(s)),-V_x(s,\bar{X}^{t,x;\bar{u}}(s)),-V_{xx}(s,\bar{X}^{t,x;\bar{u}}(s)),\bar{u}(s)\big)\\
 =&\max\limits_{u\in\mathbf{U}}G\big(s,\bar{X}^{t,x;\bar{u}}(s),-V(s,\bar{X}^{t,x;\bar{u}}(s)),
  -V_x(s,\bar{X}^{t,x;\bar{u}}(s)),-V_{xx}(s,\bar{X}^{t,x;\bar{u}}(s)),u\big),
\end{aligned}
\end{equation}
a.e.$s\in[t,T],\mathbf{P}\mbox{-}$a.s. Then
\begin{equation}\label{MP}
\begin{aligned}
&\mathcal{H}(s,\bar{X}^{t,x;\bar{u}}(s),\bar{Y}^{t,x;\bar{u}}(s),\bar{Z}^{t,x;\bar{u}}(s),\bar{u}(s),p(s),q(s),P(s))\\
=&\max\limits_{u\in\mathbf{U}}\mathcal{H}(s,\bar{X}^{t,x;\bar{u}}(s),\bar{Y}^{t,x;\bar{u}}(s),\bar{Z}^{t,x;\bar{u}}(s),u,p(s),q(s),P(s)),\ a.e. s\in[t,T],\mathbf{P}\mbox{-}a.s.
\end{aligned}
\end{equation}
\end{corollary}

{\it Proof.}\quad By the stochastic verification theorem in \cite{SY13}, from (\ref{DPP})  we know that $\bar{u}(\cdot)$ is an optimal control. Moreover, from  (\ref{generalized Hamiltonian}) and (\ref{DPP}), we have
\begin{equation*}
\begin{aligned}
    &-\big\langle V_x(s,\bar{X}^{t,x;\bar{u}}(s)),\bar{b}(s)\big\rangle
     -\frac{1}{2}\mbox{tr}\big(\bar{\sigma}(s)^\top V_{xx}(s,\bar{X}^{t,x;\bar{u}}(s))\bar{\sigma}(s)\big)\\
    &+f\big(s,\bar{X}^{t,x;\bar{u}}(s),-V(s,\bar{X}^{t,x;\bar{u}}(s)),-V_x(s,\bar{X}^{t,x;\bar{u}}(s))\bar{\sigma}(s),\bar{u}(s)\big)\\
\geq&-\big\langle V_x(s,\bar{X}^{t,x;\bar{u}}(s)),b(s,\bar{X}^{t,x;\bar{u}}(s),u)\big\rangle\\
    &-\frac{1}{2}\mbox{tr}\big(\sigma(s,\bar{X}^{t,x;\bar{u}}(s),u)^\top V_{xx}(s,\bar{X}^{t,x;\bar{u}}(s))\sigma(s,\bar{X}^{t,x;\bar{u}}(s),u)\big)\\
    &+f\big(s,\bar{X}^{t,x;\bar{u}}(s),-V(s,\bar{X}^{t,x;\bar{u}}(s)),-V_x(s,\bar{X}^{t,x;\bar{u}}(s))\sigma(s,\bar{X}^{t,x;\bar{u}}(s),u),u\big).
\end{aligned}
\end{equation*}
Noting that $p(s)=-V_x(s,\bar{X}^{t,x;\bar{u}}(s))$, thus for all $u\in\mathbf{U},\ a.e.s\in[t,T],\ \mathbf{P}\mbox{-}a.s.$,
\begin{equation*}
\begin{aligned}
    &\big\langle p(s),\bar{b}(s)\big\rangle
     +\frac{1}{2}\big\langle P(s)\bar{\sigma}(s),\bar{\sigma}(s)\big\rangle
    +f\big(s,\bar{X}^{t,x;\bar{u}}(s),-V(s,\bar{X}^{t,x;\bar{u}}(s)),\bar{\sigma}(s)^\top p(s),\bar{u}(s)\big)\\
    &-\frac{1}{2}\big\langle\big(P(s)+V_{xx}(s,\bar{X}^{t,x;\bar{u}}(s))\big)\bar{\sigma}(s),\bar{\sigma}(s)\big\rangle\\
\geq&\ \big\langle p(s),b(s,\bar{X}^{t,x;\bar{u}}(s),u)\big\rangle
     +\frac{1}{2}\big\langle P(s)\sigma(s,\bar{X}^{t,x;\bar{u}}(s),u),\sigma(s,\bar{X}^{t,x;\bar{u}}(s),u)\big\rangle\\
    &+f\big(s,\bar{X}^{t,x;\bar{u}}(s),-V(s,\bar{X}^{t,x;\bar{u}}(s)),-V_x(s,\bar{X}^{t,x;\bar{u}}(s))\sigma(s,\bar{X}^{t,x;\bar{u}}(s),u),u\big)\\
    &-\frac{1}{2}\big\langle\big(P(s)+V_{xx}(s,\bar{X}^{t,x;\bar{u}}(s))\big)
     \sigma(s,\bar{X}^{t,x;\bar{u}}(s),u),\sigma(s,\bar{X}^{t,x;\bar{u}}(s),u)\big\rangle,
\end{aligned}
\end{equation*}
where $P(\cdot)$ satisfies (\ref{adjoint equation-second order}). From the definition of $\mathcal{H}_1$ in (\ref{mathcal H1}), we get for all $u\in\mathbf{U},\ a.e.s\in[t,T],\ \mathbf{P}\mbox{-}a.s.$,
\begin{equation*}
\begin{aligned}
    &\mathcal{H}_1(s,\bar{X}^{t,x;\bar{u}}(s),\bar{u}(s))-\big\langle q(s)-P(s)\bar{\sigma}(s),\bar{\sigma}(s)\big\rangle\\
    &-\frac{1}{2}\big\langle\big(P(s)+V_{xx}(s,\bar{X}^{t,x;\bar{u}}(s))\big)\bar{\sigma}(s),\bar{\sigma}(s)\big\rangle\\
\geq&\ \mathcal{H}_1(s,\bar{X}^{t,x;\bar{u}}(s),u)-\big\langle q(s)-P(s)\bar{\sigma}(s),\sigma(s,\bar{X}^{t,x;\bar{u}}(s),u)\big\rangle)\\
    &-\frac{1}{2}\big\langle\big(P(s)+V_{xx}(s,\bar{X}^{t,x;\bar{u}}(s))\big)
     \sigma(s,\bar{X}^{t,x;\bar{u}}(s),u),\sigma(s,\bar{X}^{t,x;\bar{u}}(s),u)\big\rangle.
\end{aligned}
\end{equation*}
Then, noting that $q(s)=-V_{xx}(s,\bar{X}^{t,x;\bar{u}}(s))\bar{\sigma}(s)$, we have
\begin{equation*}
\begin{aligned}
    &\mathcal{H}_1(s,\bar{X}^{t,x;\bar{u}}(s),\bar{u}(s))-\mathcal{H}_1(s,\bar{X}^{t,x;\bar{u}}(s),u)\\
\geq&
-\frac{1}{2}\big\langle\big(P(s)+V_{xx}(s,\bar{X}^{t,x;\bar{u}}(s))\big)\big(\bar{\sigma}(s)-
     \sigma(s,\bar{X}^{t,x;\bar{u}}(s),u)\big),\bar{\sigma}(s)-\sigma(s,\bar{X}^{t,x;\bar{u}}(s),u)\big\rangle\\
\geq&\ 0,\quad\text{ for all } u\in\mathbf{U},\ a.e.s\in[t,T],\ \mathbf{P}\mbox{-}a.s.
\end{aligned}
\end{equation*}
The last inequality holds due to (\ref{relation Hu and Peng-P,Q}). Noticing that from (\ref{relation Hu and Peng-p,q}) and (\ref{1 equation}) we have $\bar{Z}^{t,x,\bar{u}}(s)=\sigma(s,x,\bar{u})p(s)$, thus (\ref{MP}) is valid by applying that $\bar{u}(\cdot)$ achieves the same maximum value in $\mathcal{H}_1$ and $\mathcal{H}$. The proof is complete.\quad$\Box$

\vspace{1mm}

{\bf Remark 3.8}\ From the above proof, one can check that if $P(s)=-V_{xx}(s,\bar{X}^{t,x;\bar{u}}(s))$, from (\ref{MP}) we can obtain (\ref{DPP}). This means that the DPP can be derived by the MP. Moreover, noticing that (\ref{DPP}) is a necessary and sufficient condition for the optimal control, thus the maximum principle (\ref{MP}) together with $P(s)=-V_{xx}(s,\bar{X}^{t,x;\bar{u}}(s))$ is also a sufficient condition for the optimal control. We mention that the maximum principle (\ref{MP}) is a necessary condition for the optimal control, but one can give lots of examples to show that $P(s)=-V_{xx}(s,\bar{X}^{t,x;\bar{u}}(s))$ is not a necessary condition for the optimal control. Finally, let us stress that $P(s)=-V_{xx}(s,\bar{X}^{t,x;\bar{u}}(s))$ plays the similar role as the additional convex/concave conditions such that maximum principle is a sufficient condition for the optimal control, see \cite{SY13,YZ99}.


\section{Proof of the Main Results}

\subsection{Proof of Theorem 3.1}

We split the proof into several steps.

\vspace{1mm}

{\it Step 1.}\ Variational equation for SDE.

Fix an $s\in[t,T]$. For any $x^1\in\mathbf{R}^n$, denote by $X^{s,x^1;\bar{u}}(\cdot)$ the solution to the following SDE on $[s,T]$:
\begin{equation}\label{disturbed controlled SDE}
{\small X^{s,x^1;\bar{u}}(r)=x^1+\int_s^rb(\alpha,X^{s,x^1;\bar{u}}(\alpha),\bar{u}(\alpha))d\alpha
 +\int_s^r\sigma(\alpha,X^{s,x^1;\bar{u}}(\alpha),\bar{u}(\alpha))dW(\alpha).}
\end{equation}
It is clear that (\ref{disturbed controlled SDE}) can be regarded as an SDE on $\big(\Omega,\mathcal{F},\{\mathcal{F}_r^t\}_{r\geq t},\mathbf{P}(\cdot|\mathcal{F}_s^t)(\omega)\big)$ for $\mathbf{P}\mbox{-}a.s. \omega$, where
$\mathbf{P}(\cdot|\mathcal{F}_s^t)(\omega)$ is the regular conditional probability given $\mathcal{F}_s^t$ defined on $(\Omega,\mathcal{F})$. For any $s\leq r\leq T$, set $\hat{X}(r):=X^{s,x^1;\bar{u}}(r)-\bar{X}^{t,x;\bar{u}}(r)$. Thus by a standard argument, we have for any integer $k\geq1$,
\begin{equation}\label{estimate of SDE}
\mathbb{E}\Big[\sup\limits_{s\leq r\leq T}\big|\hat{X}(r)\big|^{2k}\Big|\mathcal{F}_s^t\Big]\leq C|x^1-\bar{X}^{t,x;\bar{u}}(s)|^{2k},\ \mathbf{P}\mbox{-}a.s.
\end{equation}

Now we write the equation for $\hat{X}(\cdot)$ as
\begin{equation}\label{first-order variational equation-SDE}
\left\{
\begin{aligned}
d\hat{X}(r)&=\big[\bar{b}_x(r)\hat{X}(r)+\varepsilon_1(r)\big]dr+\big[\bar{\sigma}_x(r)\hat{X}(r)+\varepsilon_2(r)\big]dW(r),r\in[s,T],\\
 \hat{X}(s)&=x^1-\bar{X}^{t,x;\bar{u}}(s),
\end{aligned}
\right.
\end{equation}
and
\begin{equation}\label{second-order variational equation-SDE}
\left\{
\begin{aligned}
d\hat{X}(r)&=\big[\bar{b}_x(r)\hat{X}(r)+\frac{1}{2}\hat{X}(r)^\top\bar{b}_{xx}(r)\hat{X}(r)+\varepsilon_3(r)\big]dr\\
           &\quad+\big[\bar{\sigma}_x(r)\hat{X}(r)+\frac{1}{2}\hat{X}(r)^\top\bar{\sigma}_{xx}(r)\hat{X}(r)+\varepsilon_4(r)\big]dW(r),r\in[s,T],\\
 \hat{X}(s)&=x^1-\bar{X}^{t,x;\bar{u}}(s),
\end{aligned}
\right.
\end{equation}
where
\begin{equation*}
\left\{
\begin{aligned}
\varepsilon_1(r)&:=\int_0^1\big[b_x(r,\bar{X}^{t,x;\bar{u}}(r)+\theta\hat{X}(r),\bar{u}(r))-\bar{b}_x(r)\big]\hat{X}(r)d\theta,\\
\varepsilon_2(r)&:=\int_0^1\big[\sigma_x(r,\bar{X}^{t,x;\bar{u}}(r)+\theta\hat{X}(r),\bar{u}(r))-\bar{\sigma}_x(r)\big]\hat{X}(r)d\theta,\\
\varepsilon_3(r)&:=\int_0^1(1-\theta)\hat{X}(r)^\top\big[b_{xx}(r,\bar{X}^{t,x;\bar{u}}(r)+\theta\hat{X}(r),\bar{u}(r))-\bar{b}_{xx}(r)\big]\hat{X}(r)d\theta,\\
\varepsilon_4(r)&:=\int_0^1(1-\theta)\hat{X}(r)^\top\big[\sigma_{xx}(r,\bar{X}^{t,x;\bar{u}}(r)+\theta\hat{X}(r),\bar{u}(r))
                 -\bar{\sigma}_{xx}(r)\big]\hat{X}(r)d\theta.
\end{aligned}
\right.
\end{equation*}
and $\hat{X}(r)^\top\bar{\varphi}_{xx}(r)\hat{X}(r)\equiv\big(\hat{X}(r)^\top\bar{\varphi}_{xx}^1(r)\hat{X}(r),\cdots,\hat{X}(r)^\top\bar{\varphi}_{xx}^n(r)\hat{X}(r)\big)^\top$ for $\varphi=b,\sigma$.

\vspace{1mm}

{\it Step 2.}\quad Estimates of remainder terms of SDE.

For any integer $k\geq1$, we have
\begin{equation}\label{Estimates of remainder terms of SDE}
\left\{
\begin{aligned}
&\mathbb{E}\Big[\int_s^T|\varepsilon_1(r)|^{2k}dr\big|\mathcal{F}_s^t\Big]=o(|x^1-\bar{X}^{t,x;\bar{u}}(s)|^{2k}),\ \mathbf{P}\mbox{-}a.s.,\\
&\mathbb{E}\Big[\int_s^T|\varepsilon_2(r)|^{2k}dr\big|\mathcal{F}_s^t\Big]=o(|x^1-\bar{X}^{t,x;\bar{u}}(s)|^{2k}),\ \mathbf{P}\mbox{-}a.s.,\\
&\mathbb{E}\Big[\int_s^T|\varepsilon_3(r)|^kdr\big|\mathcal{F}_s^t\Big]=o(|x^1-\bar{X}^{t,x;\bar{u}}(s)|^{2k}),\ \mathbf{P}\mbox{-}a.s.,\\
&\mathbb{E}\Big[\int_s^T|\varepsilon_4(r)|^kdr\big|\mathcal{F}_s^t\Big]=o(|x^1-\bar{X}^{t,x;\bar{u}}(s)|^{2k}),\ \mathbf{P}\mbox{-}a.s.
\end{aligned}
\right.
\end{equation}
Here $\mathbb{E}\Big[\int_s^T|\varepsilon_1(r)|^{2k}dr\big|\mathcal{F}_s^t\Big]=o(|x^1-\bar{X}^{t,x;\bar{u}}(s)|^{2k}),\ \mathbf{P}\mbox{-}a.s.,$ means that for $\mathbf{P}\mbox{-}a.s.$ $\omega$ fixed, $\mathbb{E}\Big[\int_s^T|\varepsilon_1(r)|^{2k}dr\big|\mathcal{F}_s^t\Big](\omega)=
o(|x^1-\bar{X}^{t,x;\bar{u}}(s,\omega)|^{2k})$, where $o(\cdot)$ is almost surely a deterministic function under the regular conditional probability $\mathbf{P}(\cdot|\mathcal{F}_s^t)(\omega)$. Moreover, $o(|x^1-\bar{X}^{t,x;\bar{u}}(s,\omega)|^{2k})$ depends only on the size of $|x^1-\bar{X}^{t,x;\bar{u}}(s,\omega)|$, and it is independent of $x^1$. Such notation has similar meaning for other estimates in (\ref{Estimates of remainder terms of SDE}) as well as in the sequel of the paper.

To prove (\ref{Estimates of remainder terms of SDE}), by the continuity and uniformly boundedness of $b_x,b_{xx},\sigma_x,\sigma_{xx}$ and dominated convergence theorem as well as (\ref{estimate of SDE}), we have
\begin{equation*}
{\small\begin{aligned}
&\mathbb{E}\Big[\int_s^T|\varepsilon_1(r)|^{2k}dr\big|\mathcal{F}_s^t\Big]
\leq\int_s^T\Big\{\mathbb{E}\Big[\int_0^1\big|b_x(r,\bar{X}^{t,x;\bar{u}}(r)+\theta\hat{X}(r),\bar{u}(r))
 -\bar{b}_x(r)\big|^{4k}d\theta\Big|\mathcal{F}_s^t\Big]\Big\}^{\frac{1}{2}}\\
&\qquad\cdot\Big\{\mathbb{E}\Big[\big|\hat{X}(r)\big|^{4k}\Big|\mathcal{F}_s^t\Big]\Big\}^{\frac{1}{2}}dr=o(|x^1-\bar{X}^{t,x;\bar{u}}(s)|^{2k}),\ \mathbf{P}\mbox{-}a.s.
\end{aligned}}
\end{equation*}
Thus the first equality in (\ref{Estimates of remainder terms of SDE}) holds, and similar for the second one. Moreover, from the modulus continuity of $b_{xx}$, we have
\begin{equation*}
{\small\begin{aligned}
&\mathbb{E}\Big[\int_s^T|\varepsilon_3(r)|^kdr\big|\mathcal{F}_s^t\Big]
\leq\int_s^T\Big\{\mathbb{E}\Big[\int_0^1\big|b_{xx}(r,\bar{X}^{t,x;\bar{u}}(r)+\theta\hat{X}(r),\bar{u}(r))
 -\bar{b}_{xx}(r)\big|^{2k}d\theta\Big|\mathcal{F}_s^t\Big]\Big\}^{\frac{1}{2}}\\
&\qquad\cdot\Big\{\mathbb{E}\Big[\big|\hat{X}(r)\big|^{4k}\Big|\mathcal{F}_s^t\Big]\Big\}^{\frac{1}{2}}dr=o(|x^1-\bar{X}^{t,x;\bar{u}}(s)|^{2k}),\ \mathbf{P}\mbox{-}a.s.,
\end{aligned}}
\end{equation*}
the third equality in (\ref{Estimates of remainder terms of SDE}) holds, and the fourth one can be proved similarly.

\vspace{1mm}

{\it Step 3.}\quad Duality relation.

Applying It\^{o}'s formula to $\big\langle p(\cdot),\hat{X}(\cdot)\big\rangle$, by (\ref{adjoint equation-first order}), (\ref{second-order variational equation-SDE}), we have
\begin{equation}\label{applying Ito's formula}
\begin{aligned}
&d\big\langle p(r),\hat{X}(r)\big\rangle=\big[-\big\langle\hat{X}(r),\bar{f}_y(r)p(r)+\bar{f}_z(r)\bar{\sigma}_x(r)^\top p(r)+\bar{f}_z(r)q(r)+\bar{f}_x(r)\big\rangle\\
&\quad+\frac{1}{2}\big\langle p(r),\hat{X}(r)^\top\bar{b}_{xx}(r)\hat{X}(r)\big\rangle+\frac{1}{2}\big\langle q(r),\hat{X}(r)^\top\bar{\sigma}_{xx}(r)\hat{X}(r)\big\rangle\\
&\quad+\big\langle p(r),\varepsilon_3(r)\big\rangle+\big\langle q(r),\varepsilon_4(r)\big\rangle\big]dr+\big[\big\langle q(r),\hat{X}(r)\big\rangle+\big\langle p(r),\bar{\sigma}_x(r)\hat{X}(r)\\
&\quad+\frac{1}{2}\hat{X}(r)^\top\bar{\sigma}_{xx}(r)\hat{X}(r)\big\rangle+\big\langle p(r),\varepsilon_4(r)\big\rangle\big]dW(r),\ r\in[s,T].
\end{aligned}
\end{equation}
Setting $A(r):=\hat{X}(r)\hat{X}(r)^\top$ and applying It\^{o}'s formula to $\mbox{tr}\big\{P(\cdot)A(\cdot)\big\}$, using (\ref{adjoint equation-second order}), we obtain
\begin{equation}\label{applying Ito's formula-3}
{\small\begin{aligned}
&d\mbox{tr}\big\{P(r)A(r)\big\}
=\mbox{tr}\Big\{-A(r)\bar{f}_y(r)P(r)-A(r)\bar{f}_z(r)\bar{\sigma}_x(r)^\top P(r)-A(r)P(r)\bar{f}_z(r)\bar{\sigma}_x(r)\\
&-A(r)\bar{f}_z(r)Q(r)-A(r)\bar{b}_{xx}(r)^\top p(r)-A(r)\bar{\sigma}_{xx}(r)^\top\big[\bar{f}_z(r)p(r)+q(r)\big]+P(r)\varepsilon_5(r)-A(r)\\
&\times\big[I_{n\times n},p(r),\bar{\sigma}_x(r)^\top p(r)+q(r)\big]D^2\bar{f}(r)\big[I_{n\times n},p(r),\bar{\sigma}_x(r)^\top p(r)+q(r)\big]^\top+Q(r)\varepsilon_6(r)\Big\}dr\\
&+\mbox{tr}\Big\{A(r)Q(r)+P(r)A(r)\bar{\sigma}_x(r)^\top+P(r)\bar{\sigma}_x(r)A(r)+P(r)\varepsilon_6(r)\Big\}dW(r),\ r\in[s,T].
\end{aligned}}
\end{equation}
where
\begin{equation*}
\left\{
{\small\begin{aligned}
\varepsilon_5(r)&:=\hat{X}(r)\varepsilon_1(r)^\top+\varepsilon_1(r)\hat{X}(r)^\top+\bar{\sigma}_x(r)\hat{X}(r)\varepsilon_2(r)^\top
                +\varepsilon_2(r)\hat{X}(r)^\top\bar{\sigma}_x(r)^\top+\varepsilon_2(r)\varepsilon_2(r)^\top,\\
\varepsilon_6(r)&:=\hat{X}(r)\varepsilon_2(r)^\top+\varepsilon_2(r)\hat{X}(r)^\top.
\end{aligned}}
\right.
\end{equation*}
By (\ref{applying Ito's formula}) and (\ref{applying Ito's formula-3}), for $\hat{Y}(r):=\big\langle p(r),\hat{X}(r)\big\rangle+\frac{1}{2}\big\langle P(r)\hat{X}(r),\hat{X}(r)\big\rangle$, we obtain
\begin{equation}\label{applying Ito's formula-4}
\begin{aligned}
d\hat{Y}(r)=C(r)dr+\hat{Z}(r)dW(r),\ r\in[s,T],
\end{aligned}
\end{equation}
where
\begin{equation*}
{\small\begin{aligned}
\hat{Z}(r):=&\ \big\langle q(r),\hat{X}(r)\big\rangle+\big\langle p(r),\bar{\sigma}_x(r)^\top\hat{X}(r)+\frac{1}{2}\hat{X}(r)^\top\bar{\sigma}_{xx}(r)\hat{X}(r)\big\rangle+\big\langle p(r),\varepsilon_4(r)\big\rangle\\
            &+\frac{1}{2}\mbox{tr}\Big\{A(r)Q(r)+P(r)A(r)\bar{\sigma}_x(r)^\top+P(r)\bar{\sigma}_x(r)A(r)+P(r)\varepsilon_6(r)\Big\},\\
      C(r):=&-\big\langle\hat{X}(r),\bar{f}_y(r)p(r)+\bar{f}_z(r)\bar{\sigma}_x(r)^\top p(r)+\bar{f}_z(r)q(r)+\bar{f}_x(r)\big\rangle+\big\langle p(r),\varepsilon_3(r)\big\rangle\\
            &+\big\langle q(r),\varepsilon_4(r)\big\rangle+\mbox{tr}\Big\{-\frac{1}{2}A(r)\bar{f}_y(r)P(r)-A(r)\bar{f}_z(r)\bar{\sigma}_x(r)^\top P(r)\\
            &-\frac{1}{2}A(r)\bar{f}_z(r)Q(r)-\frac{1}{2}A(r)\bar{\sigma}_{xx}(r)^\top\bar{f}_z(r)p(r)+\frac{1}{2}P(r)\varepsilon_5(r)+\frac{1}{2}Q(r)\varepsilon_6(r)\\
            &-\frac{1}{2}A(r)\big[I_{n\times n},p(r),\bar{\sigma}_x(r)^\top p(r)+q(r)\big]D^2\bar{f}(r)\big[I_{n\times n},p(r),\bar{\sigma}_x(r)^\top p(r)+q(r)\big]^\top\Big\}.
\end{aligned}}
\end{equation*}

\vspace{1mm}

{\it Step 4.}\quad Variational equation for BSDE.

For the above $x^1\in\mathbf{R}^n$, recall $X^{s,x^1;\bar{u}}(\cdot)$ is given by (\ref{disturbed controlled SDE}) and denote by $(Y^{s,x^1;\bar{u}}(\cdot),\\Z^{s,x^1;\bar{u}}(\cdot))$ the solution to the following BSDE on $[s,T]$:
\begin{equation}\label{disturbed controlled BSDE}
\begin{aligned}
 Y^{s,x^1;\bar{u}}(r)&=\phi(X^{s,x^1;\bar{u}}(T))
                      +\int_r^Tf(\alpha,X^{s,x^1;\bar{u}}(\alpha),Y^{s,x^1;\bar{u}}\alpha),Z^{s,x^1;\bar{u}}(\alpha),\bar{u}(\alpha))d\alpha\\
                     &\quad-\int_r^TZ^{s,x^1;\bar{u}}(\alpha)dW(\alpha),
\end{aligned}
\end{equation}
and similarly (\ref{disturbed controlled BSDE}) is a BSDE on $\big(\Omega,\mathcal{F},\{\mathcal{F}_r^t\}_{r\geq t},\mathbf{P}(\cdot|\mathcal{F}_s^t)(\omega)\big)$ for $\mathbf{P}\mbox{-}a.s. \omega$.

For any $s\leq r\leq T$, set
\begin{equation}\label{widetilde Y, widetilde Z}
\widetilde{Y}(r):=Y^{s,x^1;\bar{u}}(r)-\hat{Y}(r),\ \widetilde{Z}(r):=Z^{s,x^1;\bar{u}}(r)-\hat{Z}(r).
\end{equation}
Thus by (\ref{applying Ito's formula-4}) and (\ref{disturbed controlled BSDE}), we get
%
\begin{equation}\label{minus}
\left\{
{\small\begin{aligned}
&d\big(\widetilde{Y}(r)-\bar{Y}^{t,x;\bar{u}}(r)\big)=-\big[C(r)+f(r,X^{s,x^1;\bar{u}}(r),Y^{s,x^1;\bar{u}}(r),Z^{s,x^1;\bar{u}}(r),\bar{u}(r))\\
&\quad-f(r,\bar{X}^{t,x;\bar{u}}(r),\bar{Y}^{t,x;\bar{u}}(r),\bar{Z}^{t,x;\bar{u}}(r),\bar{u}(r))\big]dr+\big[\widetilde{Z}(r)-\bar{Z}^{t,x;\bar{u}}(r)\big]dW(r),\\
&\widetilde{Y}(T)-\bar{Y}^{t,x;\bar{u}}(T)=\int_0^1(1-\theta)\hat{X}(T)^\top\big[\phi_{xx}(\bar{X}^{t,x;\bar{u}}(T)+\theta\hat{X}(T))\\
&\quad-\phi_{xx}(\bar{X}^{t,x;\bar{u}}(T))\big]\hat{X}(T)d\theta.
\end{aligned}}
\right.
\end{equation}
By the modulus continuity of $\phi_{xx}$, we have
\begin{equation*}
\begin{aligned}
&\mathbb{E}\big[\big|\widetilde{Y}(T)-\bar{Y}^{t,x;\bar{u}}(T)\big|\big|\mathcal{F}_s^t\big]
\leq\Big\{\mathbb{E}\Big[\Big(\int_0^1\big|\phi_{xx}(\bar{X}^{t,x;\bar{u}}(T)+\theta\hat{X}(T))\\
&-\phi_{xx}(\bar{X}^{t,x;\bar{u}}(T))\big|d\theta\Big)^2\Big|\mathcal{F}_s^t\Big]\Big\}^{\frac{1}{2}}\Big\{\mathbb{E}\Big[|\hat{X}(T)|^4\Big|\mathcal{F}_s^t\Big]\Big\}^{\frac{1}{2}}=o(|x^1-\bar{X}^{t,x;\bar{u}}(s)|^2),\ \mathbf{P}\mbox{-}a.s.
\end{aligned}
\end{equation*}
Noting (\ref{widetilde Y, widetilde Z}), we have
\begin{equation*}
\begin{aligned}
&f(r,X^{s,x^1;\bar{u}}(r),Y^{s,x^1;\bar{u}}(r),Z^{s,x^1;\bar{u}}(r),\bar{u}(r))
-f(r,\bar{X}^{t,x;\bar{u}}(r),\bar{Y}^{t,x;\bar{u}}(r),\bar{Z}^{t,x;\bar{u}}(r),\bar{u}(r))\\
&=f(r,\bar{X}^{s,x^1;\bar{u}}(r)+\hat{X}(r),\widetilde{Y}(r)+\hat{Y}(r),\widetilde{Z}(r)+\hat{Z}(r),\bar{u}(r))\\
&\quad-f(r,\bar{X}^{s,x^1;\bar{u}}(r)+\hat{X}(r),\bar{Y}^{t,x;\bar{u}}(r)+\hat{Y}(r),\bar{Z}^{t,x;\bar{u}}(r)+\hat{Z}(r),\bar{u}(r))\\
&\quad+f(r,\bar{X}^{s,x^1;\bar{u}}(r)+\hat{X}(r),\bar{Y}^{t,x;\bar{u}}(r)+\hat{Y}(r),\bar{Z}^{t,x;\bar{u}}(r)+\hat{Z}(r),\bar{u}(r))\\
&\quad-f(r,\bar{X}^{t,x;\bar{u}}(r),\bar{Y}^{t,x;\bar{u}}(r),\bar{Z}^{t,x;\bar{u}}(r),\bar{u}(r))\\
&=\widetilde{f}_y(r)\big(\widetilde{Y}(r)-\bar{Y}^{t,x;\bar{u}}(r)\big)+\widetilde{f}_z(r)\big(\widetilde{Z}(r)
 -\bar{Z}^{t,x;\bar{u}}(r)\big)+\bar{f}_x(r)\hat{X}(r)\\
&\quad+\bar{f}_y(r)\hat{Y}(r)+\bar{f}_z(r)\hat{Z}(r)+[\hat{X}(r)^\top,\hat{Y}(r),\hat{Z}(r)]\widetilde{D}^2f(r)[\hat{X}(r)^\top,\hat{Y}(r),\hat{Z}(r)]^\top,
\end{aligned}
\end{equation*}
where
\begin{equation*}
\left\{
\begin{aligned}
 \widetilde{f}_y(r)&:=\int_0^1f_y\big(r,\bar{X}^{s,x^1;\bar{u}}(r)+\hat{X}(r),\bar{Y}^{t,x;\bar{u}}(r)
                    +\hat{Y}(r)+\theta(\widetilde{Y}(r)-\bar{Y}^{t,x;\bar{u}}(r)),\\
                   &\qquad\qquad\bar{Z}^{t,x;\bar{u}}(r)+\hat{Z}(r)+\theta(\widetilde{Z}(r)-\bar{Z}^{t,x;\bar{u}}(r)),\bar{u}(r)\big)d\theta,\\
 \widetilde{f}_z(r)&:=\int_0^1f_z\big(r,\bar{X}^{s,x^1;\bar{u}}(r)+\hat{X}(r),\bar{Y}^{t,x;\bar{u}}(r)
                    +\hat{Y}(r)+\theta(\widetilde{Y}(r)-\bar{Y}^{t,x;\bar{u}}(r)),\\
                   &\qquad\qquad\bar{Z}^{t,x;\bar{u}}(r)+\hat{Z}(r)+\theta(\widetilde{Z}(r)-\bar{Z}^{t,x;\bar{u}}(r)),\bar{u}(r)\big)d\theta,\\
\widetilde{D}^2f(r)&:=\int_0^1\int_0^1\lambda D^2f\big(r,\bar{X}^{t,x;\bar{u}}(r)+\lambda\theta\hat{X}(r),\bar{Y}^{t,x;\bar{u}}(r)+\lambda\theta\hat{Y}(r),\\
                   &\qquad\qquad\bar{Z}^{t,x;\bar{u}}(r)+\lambda\theta\hat{Z}(r),\bar{u}(r)\big)d\lambda d\theta.
\end{aligned}
\right.
\end{equation*}
Using the definition of  $\hat{Y}(r)$ and $\hat{Z}(r)$, we obtain
\begin{equation}\label{minus-integral form}
{\small\begin{aligned}
&\widetilde{Y}(s)-\bar{Y}^{t,x;\bar{u}}(s)
=o(|x^1-\bar{X}^{t,x;\bar{u}}(s)|^2)+\int_s^T\bigg\{\widetilde{f}_y(r)\big(\widetilde{Y}(r)-\bar{Y}^{t,x;\bar{u}}(r)\big)\\
&+\widetilde{f}_z(r)\big(\widetilde{Z}(r)-\bar{Z}^{t,x;\bar{u}}(r)\big)+[\hat{X}(r)^\top,\hat{Y}(r),\hat{Z}(r)]\widetilde{D}^2f(r)[\hat{X}(r)^\top,\hat{Y}(r),\hat{Z}(r)]^\top\\
&-\frac{1}{2}\mbox{tr}\Big\{\hat{X}(r)\hat{X}(r)^\top\big[I_{n\times n},p(r),\bar{\sigma}_x^\top(r)p(r)+q(r)\big]D^2\bar{f}(r)\big[I_{n\times n},p(r),\bar{\sigma}_x^\top(r)p(r)+q(r)\big]^\top\\
&+C_1(r)\Big\}\bigg\}dr-\int_s^T\big[\widetilde{Z}(r)-\bar{Z}^{t,x;\bar{u}}(r)\big]dW(r),\ \mathbf{P}\mbox{-}a.s.,
\end{aligned}}
\end{equation}
where
\begin{equation*}
\begin{aligned}
C_1(r)&:=\big\langle p(r),\varepsilon_3(r)\big\rangle+\big\langle q(r),\varepsilon_4(r)\big\rangle+\bar{f}_z(r)\big\langle p(r),\varepsilon_4(r)\big\rangle\\
      &\quad+\frac{1}{2}\mbox{tr}\big\{P(r)\varepsilon_5(r)+Q(r)\varepsilon_6(r)+\bar{f}_z(r)P(r)\varepsilon_6(r)\big\}.
\end{aligned}
\end{equation*}
Denoting the $n+2$-dimensional random vectors
\begin{equation*}
M(r):=[\hat{X}(r)^\top,\hat{Y}(r),\hat{Z}(r)],\ N(r):=\big[\hat{X}(r)^\top,\hat{X}(r)^\top p(r),\hat{X}(r)^\top\big(\bar{\sigma}_x(r)^\top p(r)+q(r)\big)\big],
\end{equation*}
we obtain that
\begin{equation}\label{minus-integral form-simple}
\begin{aligned}
&\widetilde{Y}(s)-\bar{Y}^{t,x;\bar{u}}(s)
=o(|x^1-\bar{X}^{t,x;\bar{u}}(s)|^2)+\int_s^T\bigg\{\widetilde{f}_y(r)\big(\widetilde{Y}(r)-\bar{Y}^{t,x;\bar{u}}(r)\big)\\
&\quad+\widetilde{f}_z(r)\big(\widetilde{Z}(r)-\bar{Z}^{t,x;\bar{u}}(r)\big)+C_1(r)+M(r)\widetilde{D}^2f(r)M(r)^\top\\
&\quad-\frac{1}{2}N(r)D^2\bar{f}(r)N(r)^\top\Big\}\bigg\}dr-\int_s^T\big[\widetilde{Z}(r)-\bar{Z}^{t,x;\bar{u}}(r)\big]dW(r),\ \mathbf{P}\mbox{-}a.s.
\end{aligned}
\end{equation}

\vspace{1mm}

{\it Step 5.}\quad Estimates of remainder terms of BSDEs.

Noting that
\begin{equation*}
\begin{aligned}
&M(r)\widetilde{D}^2f(r)M(r)^\top-\frac{1}{2}N(r)D^2\bar{f}(r)N(r)^\top=M(r)\widetilde{D}^2f(r)M(r)^\top\\
&\quad-N(r)\widetilde{D}^2f(r)N(r)^\top+N(r)\widetilde{D}^2f(r)N(r)^\top-\frac{1}{2}N(r)D^2\bar{f}(r)N(r)^\top\\
&:=\Pi(r)\widetilde{D}^2f(r)\Pi(r)+N(r)\Big[\widetilde{D}^2f(r)-\frac{1}{2}D^2\bar{f}(r)\Big]N(r)^\top,
\end{aligned}
\end{equation*}
where $\Pi(r):=M(r)-N(r)$. We have
\begin{equation}\label{estimate Pi,N,C1}
\left\{
\begin{aligned}
&\mathbb{E}\Big[\Big(\int_s^T\big|C_1(r)\big|dr\Big)^2\Big|\mathcal{F}_s^t\Big]=o(|x^1-\bar{X}^{t,x;\bar{u}}(s)|^4),\ \mathbf{P}\mbox{-}a.s.,\\
&\mathbb{E}\Big[\Big(\int_s^T\big|\Pi(r)\widetilde{D}^2f(r)\Pi(r)\big|dr\Big)^2\Big|\mathcal{F}_s^t\Big]=o(|x^1-\bar{X}^{t,x;\bar{u}}(s)|^4),\ \mathbf{P}\mbox{-}a.s.,\\
&\mathbb{E}\Big[\Big(\int_s^T\Big|N(r)\Big[\widetilde{D}^2f(r)-\frac{1}{2}D^2\bar{f}(r)\Big]N(r)^\top\Big|dr\Big)^2\Big|\mathcal{F}_s^t\Big]=o(|x^1-\bar{X}^{t,x;\bar{u}}(s)|^4),\ \mathbf{P}\mbox{-}a.s.
\end{aligned}
\right.
\end{equation}
Indeed, by the boundedness of $f_z$, we have
\begin{equation*}
\begin{aligned}
&\mathbb{E}\Big[\Big(\int_s^T\big|C_1(r)\big|dr\Big)^2\Big|\mathcal{F}_s^t\Big]
\leq C\mathbb{E}\Big[\Big(\int_s^T\big\langle p(r),\varepsilon_3(r)\big\rangle dr\Big)^2\Big|\mathcal{F}_s^t\Big]\\
&\quad+C\mathbb{E}\Big[\Big(\int_s^T\big\langle q(r),\varepsilon_4(r)\big\rangle dr\Big)^2\Big|\mathcal{F}_s^t\Big]\
 +C\mathbb{E}\Big[\Big(\int_s^T\big\langle p(r),\varepsilon_4(r)\big\rangle dr\Big)^2\Big|\mathcal{F}_s^t\Big]\\
\end{aligned}
\end{equation*}
\begin{equation*}
\begin{aligned}&\quad+C\mathbb{E}\Big[\Big(\int_s^T\mbox{tr}\big\{P(r)\varepsilon_5(r)\big\}dr\Big)^2\Big|\mathcal{F}_s^t\Big]
 +C\mathbb{E}\Big[\Big(\int_s^T\mbox{tr}\big\{Q(r)\varepsilon_6(r)\big\}dr\Big)^2\Big|\mathcal{F}_s^t\Big]\\
&\quad+C\mathbb{E}\Big[\Big(\int_s^T\mbox{tr}\big\{\bar{f}_z(r)P(r)\varepsilon_6(r)\big\}dr\Big)^2\Big|\mathcal{F}_s^t\Big]:=I_1+I_2+I_3+I_4+I_5+I_6.
\end{aligned}
\end{equation*}
In the sequel, by (\ref{Estimates of remainder terms of SDE}), we have
\begin{equation*}
\begin{aligned}
I_1&
    \leq C\mathbb{E}\Big[\sup\limits_{s\leq r\leq T}\big|p(r)\big|^2\int_s^T|\varepsilon_3(r)|^2dr\Big|\mathcal{F}_s^t\Big]
    \leq C\mathbb{E}\Big[\int_s^T|\varepsilon_3(r)|^2dr\Big|\mathcal{F}_s^t\Big]\\
   &=o(|x^1-\bar{X}^{t,x;\bar{u}}(s)|^4),\ \mathbf{P}\mbox{-}a.s.,
\end{aligned}
\end{equation*}
and similar for $I_3$; by H\"{o}lder's inequality, we get
\begin{equation*}
\begin{aligned}
I_2&
    \leq C\mathbb{E}\Big[\int_s^T\big|q(r)\big|^2dr\int_s^T|\varepsilon_4(r)|^2dr\Big|\mathcal{F}_s^t\Big]\leq C\Big\{\mathbb{E}\Big[\Big(\int_s^T\big|q(r)\big|^2dr\Big)^2\Big|\mathcal{F}_s^t\Big]\Big\}^{\frac{1}{2}}\\
   &\quad\cdot\Big\{\mathbb{E}\Big[\int_s^T|\varepsilon_3(r)|^4dr\Big|\mathcal{F}_s^t\Big]\Big\}^{\frac{1}{2}}=o(|x^1-\bar{X}^{t,x;\bar{u}}(s)|^4),\ \mathbf{P}\mbox{-}a.s.;
\end{aligned}
\end{equation*}
from the definition of $\varepsilon_5(r)$ and (\ref{estimate of SDE}), (\ref{Estimates of remainder terms of SDE}), we obtain
\begin{equation*}
\begin{aligned}
I_4
   &\leq C\mathbb{E}\Big[\sup\limits_{s\leq r\leq T}\big|P(r)\hat{X}(r)\big|^2\int_s^T\big(|\varepsilon_1(r)|^2
    +|\varepsilon_2(r)|^2+|\varepsilon_2(r)|^4\big)dr\Big|\mathcal{F}_s^t\Big]\\
   &\leq C\Big\{\mathbb{E}\Big[\sup\limits_{s\leq r\leq T}\big|\hat{X}(r)\big|^8\Big|\mathcal{F}_s^t\Big]\Big\}^{1\over4}
    \Big\{\mathbb{E}\Big[\int_s^T\big(|\varepsilon_1(r)|^2+|\varepsilon_2(r)|^2+|\varepsilon_2(r)|^4\big)^4dr\Big|\mathcal{F}_s^t\Big]\Big\}^{1\over4}\\
   &=o(|x^1-\bar{X}^{t,x;\bar{u}}(s)|^4),\ \mathbf{P}\mbox{-}a.s.;
\end{aligned}
\end{equation*}
from the definition of $\varepsilon_6(r)$ and (\ref{estimate of SDE}), (\ref{Estimates of remainder terms of SDE}), we obtain
\begin{equation*}
\begin{aligned}
I_5
   &\leq C\mathbb{E}\Big[\sup\limits_{s\leq r\leq T}\big|\hat{X}(r)\big|^2\Big(\int_s^T|Q(r)\varepsilon_2(r)|dr\Big)^2\Big|\mathcal{F}_s^t\Big]\\
   &\leq C\Big\{\mathbb{E}\Big[\sup\limits_{s\leq r\leq T}\big|\hat{X}(r)\big|^4\Big|\mathcal{F}_s^t\Big]\Big\}^{1\over2}
    \Big\{\mathbb{E}\Big[\Big(\int_s^T|Q(r)\varepsilon_2(r)|dr\Big)^4\Big|\mathcal{F}_s^t\Big]\Big\}^{1\over2}\\
   &\leq C\Big\{\mathbb{E}\Big[\sup\limits_{s\leq r\leq T}\big|\hat{X}(r)\big|^4\Big|\mathcal{F}_s^t\Big]\Big\}^{1\over2}
    \Big\{\mathbb{E}\Big[\Big(\int_s^T|Q(r)|^2dr\Big)^2\Big(\int_s^T|\varepsilon_2(r)|^2dr\Big)^2\Big|\mathcal{F}_s^t\Big]\Big\}^{1\over2}\\
   &\leq C\Big\{\mathbb{E}\Big[\sup\limits_{s\leq r\leq T}\big|\hat{X}(r)\big|^4\Big|\mathcal{F}_s^t\Big]\Big\}^{1\over2}
    \Big\{\mathbb{E}\Big[\Big(\int_s^T|\varepsilon_2(r)|^8dr\Big)\Big|\mathcal{F}_s^t\Big]\Big\}^{1\over4}\\
   &=o(|x^1-\bar{X}^{t,x;\bar{u}}(s)|^4),\ \mathbf{P}\mbox{-}a.s.;
\end{aligned}
\end{equation*}
finally the inequality for $I_6$ can be proved similarly.
Combining the above estimates, we obtain that the first equality of (\ref{estimate Pi,N,C1}) holds.

For the second one, since
\begin{equation*}
\begin{aligned}
\Pi(r)&=M(r)-N(r)\\
      &=\big[O_{n\times n},\hat{Y}(r)-\hat{X}(r)^\top p(r),\hat{Z}(r)-\hat{X}(r)^\top\big(\bar{\sigma}_x(r)^\top p(r)+q(r)\big)\big]\\
      &=\Big[O_{n\times n},\frac{1}{2}\big\langle P(r)\hat{X}(r),\hat{X}(r)\big\rangle,\frac{1}{2}\big\langle p(r),\hat{X}(r)^\top\bar{\sigma}_{xx}(r)\hat{X}(r)\big\rangle
       +\big\langle p(r),\varepsilon_4(r)\big\rangle\\
      &\quad+\mbox{tr}\Big\{\frac{1}{2}\hat{X}(r)\hat{X}(r)^\top Q(r)+P(r)\hat{X}(r)\hat{X}(r)^\top\bar{\sigma}_x(r)^\top+\frac{1}{2}P(r)\varepsilon_6(r)\Big\}\Big],
\end{aligned}
\end{equation*}
then by the definitions of $\varepsilon_4(r),\varepsilon_6(r)$, the boundedness of $\bar{\sigma}_x,\bar{\sigma}_{xx},\widetilde{D}^2f$, and the square-integrability of $P(\cdot),Q(\cdot)$, we obtain
\begin{equation*}
{\footnotesize\begin{aligned}
&\mathbb{E}\Big[\Big(\int_s^T\big|\Pi(r)\widetilde{D}^2f(r)\Pi(r)\big|dr\Big)^2\Big|\mathcal{F}_s^t\Big]\leq C\mathbb{E}\Big[\Big(\int_s^T\big|\Pi(r)^2\big|dr\Big)^2\Big|\mathcal{F}_s^t\Big]\\
&\leq C\mathbb{E}\Big[\Big(\int_s^T\Big[\big|\hat{X}(r)\big|^4\big(|P(r)|^2+|Q(r)|^2+|p(r)|^2\big)+|p(r)|^2|\varepsilon_4(r)|^2+|P(r)|^2|\varepsilon_6(r)|^2\big)\Big]dr\Big)^2\Big|\mathcal{F}_s^t\Big]\\
&\leq C\mathbb{E}\Big[\Big(\int_s^T\Big[\big|\hat{X}(r)\big|^4\big(|P(r)|^2+|Q(r)|^2+|p(r)|^2\big)\Big]dr\Big)^2\Big|\mathcal{F}_s^t\Big]+o(|x^1-\bar{X}^{t,x;\bar{u}}(s)|^4)\\
&\leq C\Big(\mathbb{E}\Big[\sup\limits_{s\leq r\leq T}\big|\hat{X}(r)\big|^{16}\Big|\mathcal{F}_s^t\Big]\Big)^{\frac{1}{2}}+o(|x^1-\bar{X}^{t,x;\bar{u}}(s)|^4)=o(|x^1-\bar{X}^{t,x;\bar{u}}(s)|^4),\ \mathbf{P}\mbox{-}a.s.
\end{aligned}}
\end{equation*}
Thus the second equality of (\ref{estimate Pi,N,C1}) follows. Finally, we prove the third one. We get
\begin{equation*}
{\footnotesize\begin{aligned}
&\mathbb{E}\Big[\Big(\int_s^T\Big|N(r)\Big[\widetilde{D}^2f(r)-\frac{1}{2}D^2\bar{f}(r)\Big]N(r)^\top\Big|dr\Big)^2
 \Big|\mathcal{F}_s^t\Big]\\
&\leq C\mathbb{E}\Big[\Big(\int_s^T\big|\hat{X}(r)\big|^2\Big(1+\big|p(r)\big|^2+\big|\bar{\sigma}_x^\top(r)p(r)+q(r)\big|^2\Big)\Big|\widetilde{D}^2f(r)-\frac{1}{2}D^2\bar{f}(r)\Big|dr\Big)^2\Big|\mathcal{F}_s^t\Big]\\
&\leq C\mathbb{E}\Big[\sup\limits_{s\leq r\leq T}\big|\hat{X}(r)\big|^4\Big(\int_s^T\big(1+\big|p(r)\big|^2+\big|q(r)\big|^2\big)
 \Big|\widetilde{D}^2f(r)-\frac{1}{2}D^2\bar{f}(r)\Big|dr\Big)^2\Big|\mathcal{F}_s^t\Big]\\
&\leq C\Big(\mathbb{E}\Big[\sup\limits_{s\leq r\leq T}\big|\hat{X}(r)\big|^8\Big|\mathcal{F}_s^t\Big]\Big)^{\frac{1}{2}}\Big(\mathbb{E}\Big[\Big(\int_s^T\big(1+\big|p(r)\big|^2+\big|q(r)\big|^2\big)\big|\widetilde{D}^2f(r)
 -\frac{1}{2}D^2\bar{f}(r)\big|dr\Big)^4\Big|\mathcal{F}_s^t\Big]\Big)^{\frac{1}{2}},\\
 &\leq C\Big(\mathbb{E}\Big[\sup\limits_{s\leq r\leq T}\big|\hat{X}(r)\big|^8\Big|\mathcal{F}_s^t\Big]\Big)^{\frac{1}{2}}\Big(\mathbb{E}\big[
 \Big(\int_s^T\big(1+\big|p(r)\big|^2+\big|q(r)\big|^2\big)^2dr\Big)^4\Big|\mathcal{F}_s^t\Big]\Big)^{\frac{1}{4}}\\
 &\qquad\cdot\Big(\mathbb{E}\Big[\int_s^T\big|\widetilde{D}^2f(r)
 -\frac{1}{2}D^2\bar{f}(r)\big|^2dr\Big)^4\Big|\mathcal{F}_s^t\Big]\Big)^{\frac{1}{4}},\ \mathbf{P}\mbox{-}a.s.
\end{aligned}}
\end{equation*}
Since $p(\cdot),q(\cdot)$ are square-integrable, by the definition of $\widetilde{D}^2f$, $\hat{Y}$ and $\hat{Z}$ and the modulus continuity of $D^2f$, we obtain the third equality of (\ref{estimate Pi,N,C1}).

By (\ref{minus-integral form-simple}), (\ref{estimate Pi,N,C1}) and Lemma 2.1, we have
\begin{equation}\label{minus-integral form-simple-high order}
{\small\begin{aligned}
&\big|\widetilde{Y}(s)-\bar{Y}^{t,x;\bar{u}}(s)\big|^2\leq C\mathbb{E}\Big[\Big(\int_s^T\big|C_1(r)\big|dr\Big)^2\Big|\mathcal{F}_s^t\Big]+C\mathbb{E}\Big[\Big(\int_s^T\big|M(r)\widetilde{D}^2f(r)M(r)^\top\\
&-\frac{1}{2}N(r)D^2f(r)N(r)^\top\big|dr\Big)^2\Big|\mathcal{F}_s^t\Big]+o(|x^1-\bar{X}^{t,x;\bar{u}}(s)|^4)\leq o(|x^1-\bar{X}^{t,x;\bar{u}}(s)|^4),\ \mathbf{P}\mbox{-}a.s.
\end{aligned}}
\end{equation}

\vspace{1mm}

{\it Step 6.}\quad Completion of the proof.

We call $x^1\in\mathbf{R}^n$ a {\it rational vector} if all its coordinate are rational numbers. Since the set of all rational vectors $x^1\in\mathbf{R}^n$ is countable, we can find a subset $\Omega_0\subseteq\Omega$ with $\mathbf{P}(\Omega_0)=1$ such that for any $\omega_0\in\Omega_0$,
\begin{equation*}
\left\{
\begin{aligned}
&V(s,\bar{X}^{t,x;\bar{u}}(s,\omega_0))=-\bar{Y}^{t,x;\bar{u}}(s,\omega_0),(\ref{estimate of SDE}),(\ref{Estimates of remainder terms of SDE}),(\ref{applying Ito's formula-4}),(\ref{minus-integral form}), (\ref{minus-integral form-simple}), (\ref{estimate Pi,N,C1}), \\
&(\ref{minus-integral form-simple-high order})\mbox{ are satisfied for any rational vector }x^1,\big(\Omega,\mathcal{F},\mathbf{P}(\cdot|\mathcal{F}_s^t)(\omega_0),W(\cdot)-W(s);\\
&u(\cdot))|_{[s,T]}\big)\in\mathcal{U}^w[s,T],\mbox{ and }\sup\limits_{s\leq r\leq T}\big[|p(r,\omega_0)|+|P(r,\omega_0)|\big]<\infty.
\end{aligned}
\right.
\end{equation*}
The first relation of the above is obtained by the DPP (see Theorem 5.4 of \cite{Peng97}). Let $\omega_0\in\Omega_0$ be fixed, then for any rational vector $x^1\in\mathbf{R}^n$, by (\ref{minus-integral form-simple-high order}), we have
\begin{equation}\label{final estimate}
\widetilde{Y}(s,\omega_0)-\bar{Y}^{t,x;\bar{u}}(s,\omega_0)=o(|x^1-\bar{X}^{t,x;\bar{u}}(s,\omega_0)|^2),\ \text{ for all } s\in[t,T].
\end{equation}
By the definition of $\widetilde{Y}(\cdot)$, we have
\begin{equation*}
\begin{aligned}
&Y^{s,x^1;\bar{u}}(s,\omega_0)-\bar{Y}^{t,x;\bar{u}}(s,\omega_0)=\big\langle p(s,\omega_0),x^1-\bar{X}^{t,x;\bar{u}}(s,\omega_0)\big\rangle\\
&+\frac{1}{2}\big\langle P(s,\omega_0)x^1-\bar{X}^{t,x;\bar{u}}(s,\omega_0),x^1-\bar{X}^{t,x;\bar{u}}(s,\omega_0)\big\rangle+o(|x^1-\bar{X}^{t,x;\bar{u}}(s,\omega_0)|^2),
\end{aligned}
\end{equation*}
for all $s\in[t,T]$. Thus
\begin{equation}\label{difference}
\begin{aligned}
    &V(s,x^1)-V(s,\bar{X}^{t,x;\bar{u}}(s,\omega_0))\leq-Y^{s,x^1;\bar{u}}(s,\omega_0)+\bar{Y}^{t,x;\bar{u}}(s,\omega_0)\\
   =&-\big\langle p(s,\omega_0),x^1-\bar{X}^{t,x;\bar{u}}(s,\omega_0)\big\rangle
     -\frac{1}{2}\big\langle P(s,\omega_0)x^1-\bar{X}^{t,x;\bar{u}}(s,\omega_0),x^1-\bar{X}^{t,x;\bar{u}}(s,\omega_0)\big\rangle\\
    &\quad+o(|x^1-\bar{X}^{t,x;\bar{u}}(s,\omega_0)|^2),\ \text{ for all } s\in[t,T].
\end{aligned}
\end{equation}
Note that the term $o(|x^1-\bar{X}^{t,x;\bar{u}}(s,\omega_0)|^2)$ in the above depends only on the size of $|x^1-\bar{X}^{t,x;\bar{u}}(s,\omega_0)|$, and it is independent of $x^1$. Therefore, by the continuity of $V(s,\cdot)$, we see that (\ref{difference}) holds for all $x^1\in\mathbf{R}^n$ (for more details see \cite{Zhou90-2}), which by definition (\ref{second-order super- and sub-jets}) proves
\begin{equation*}
\big\{-p(s)\big\}\times[-P(s),\infty\big)\in D_x^{2,+}V(s,\bar{X}^{t,x;\bar{u}}(s)),\ \text{ for all } s\in[t,T],\ \mathbf{P}\mbox{-}a.s.
\end{equation*}

Finally, fix an $\omega\in\Omega$ such that (\ref{difference}) holds for any $x^1\in\mathbf{R}^n$. For any $(\hat{p},\hat{P})\in D_x^{2,-}V(s,\bar{X}^{t,x;\bar{u}}(s))$, by definition (\ref{second-order super- and sub-jets}) we have
\begin{equation*}
\begin{aligned}
0&\leq\liminf\limits_{x^1\rightarrow\bar{X}^{t,x;\bar{u}}(s)}\left\{\frac{V(s,x^1)-V(s,\bar{X}^{t,x;\bar{u}}(s))}{|x^1-\bar{X}^{t,x;\bar{u}}(s)|}\right.\\
 &\qquad\qquad\qquad\quad\left.-\frac{\langle\hat{p},x^1-\bar{X}^{t,x;\bar{u}}(s)\rangle
  +\frac{1}{2}\big\langle \hat{P}(x^1-\bar{X}^{t,x;\bar{u}}(s)),x^1-\bar{X}^{t,x;\bar{u}}(s)\big\rangle}{|x^1-\bar{X}^{t,x;\bar{u}}(s)|}\right\}\\
 &\leq\liminf\limits_{x^1\rightarrow\bar{X}^{t,x;\bar{u}}(s)}\left\{
  -\frac{\big\langle p(s)+\hat{p},x^1-\bar{X}^{t,x;\bar{u}}(s)\big\rangle}{|x^1-\bar{X}^{t,x;\bar{u}}(s)|}\right.\\
 &\qquad\qquad\qquad\quad\left.-\frac{\frac{1}{2}\big\langle
  \big(P(s)+ \hat{P}\big)(x^1-\bar{X}^{t,x;\bar{u}}(s)),x^1-\bar{X}^{t,x;\bar{u}}(s)\big\rangle}{|x^1-\bar{X}^{t,x;\bar{u}}(s)|}\right\}\\
\end{aligned}
\end{equation*}
Then, it is necessary that
\begin{equation*}
\hat{p}=-p(s),\ \hat{P}\leq-P(s),\quad\text{ for all } s\in[t,T],\quad\mathbf{P}\mbox{-}a.s.
\end{equation*}
Thus, (\ref{connection-second order-x}) holds. The proof is complete.\quad$\Box$

\subsection{Proof of Theorem 3.3}

For any $s\in(t,T)$, take $\tau\in(s,T]$. Denote by $X^{\tau,\bar{X}^{t,x;\bar{u}}(s);\bar{u}}(\cdot)$ the solution to the following SDE on $[\tau,T]$:
\begin{equation}\label{diturbed state equation in t}
\begin{aligned}
X^{\tau,\bar{X}^{t,x;\bar{u}}(s);\bar{u}}(r)&=\bar{X}^{t,x;\bar{u}}(s)+\int_\tau^rb(\theta,X^{\tau,\bar{X}^{t,x;\bar{u}}(s);\bar{u}}(\theta),\bar{u}(\theta))
                                             d\theta\\
                                            &\quad+\int_\tau^r\sigma(\theta,X^{\tau,\bar{X}^{t,x;\bar{u}}(s);\bar{u}}(\theta),\bar{u}(\theta))dW(\theta).
\end{aligned}
\end{equation}
Set $\xi_\tau(r):=X^{\tau,\bar{X}^{t,x;\bar{u}}(s);\bar{u}}(r)-\bar{X}^{t,x;\bar{u}}(r),\tau\leq r\leq T$. We have the following estimate for any integer $k\geq1$:
\begin{equation*}
\mathbb{E}\Big[\sup\limits_{\tau\leq r\leq T}|\xi_\tau(r)|^{2k}\Big|\mathcal{F}_\tau^t\Big]\leq C\big|\bar{X}^{t,x;\bar{u}}(\tau)-\bar{X}^{t,x;\bar{u}}(s)\big|^{2k},\quad\mathbf{P}\mbox{-}a.s.
\end{equation*}
Taking $\mathbb{E}(\cdot|\mathcal{F}_s^t)$ on both sides, by a standard argument we obtain
\begin{equation}\label{estimate of SDEP in t}
\mathbb{E}\Big[\sup\limits_{\tau\leq r\leq T}|\xi_\tau(r)|^{2k}\Big|\mathcal{F}_s^t\Big]\leq C|\tau-s|^k,\quad\mathbf{P}\mbox{-}a.s.
\end{equation}
The process  $\xi_\tau(\cdot)$ satisfies the following variational equations:
\begin{equation}\label{first-order variational equation in t}
\left\{
\begin{aligned}
d\xi_\tau(r)&=\big[\bar{b}_x(r)\xi_\tau(r)+\varepsilon_{\tau1}(r)\big]dr+\big[\bar{\sigma}_x(r)\xi_\tau(r)+\varepsilon_{\tau2}(r)\big]dW(r),\ r\in[\tau,T],\\
\xi_\tau(\tau)&=-\int_s^\tau\bar{b}(r)dr-\int_s^\tau\bar{\sigma}(r)dW(r),
\end{aligned}
\right.
\end{equation}
and
\begin{equation}\label{second-order variational equation in t}
\left\{
\begin{aligned}
  d\xi_\tau(r)&=\Big\{\bar{b}_x(r)\xi_\tau(r)+\frac{1}{2}\xi_\tau(r)^\top\bar{b}_{xx}(r)\xi_\tau(r)+\varepsilon_{\tau3}(r)\Big\}dr\\
              &\qquad+\Big\{\bar{\sigma}_x(r)\xi_\tau(r)+\frac{1}{2}\xi_\tau(r)^\top\bar{\sigma}_{xx}(r)\xi_\tau(r)+\varepsilon_{\tau4}(r)\Big\}dW(r),\ r\in[\tau,T],\\
\xi_\tau(\tau)&=-\int_s^\tau\bar{b}(r)dr-\int_s^\tau\bar{\sigma}(r)dW(r),
\end{aligned}
\right.
\end{equation}
where
\begin{equation*}
\left\{
\begin{aligned}
      \varepsilon_{\tau1}(r)&:=\int_0^1\big[b_x(r,\bar{X}^{t,x;\bar{u}}(r)+\theta\xi_\tau(r),\bar{u}(r))-\bar{b}_x(r)\big]\xi_\tau(r)d\theta,\\
      \varepsilon_{\tau2}(r)&:=\int_0^1\big[\sigma_x(r,\bar{X}^{t,x;\bar{u}}(r)+\theta\xi_\tau(r),\bar{u}(r))-\bar{\sigma}_x(r)\big]\xi_\tau(r)d\theta,\\
      \varepsilon_{\tau3}(r)&:=\int_0^1(1-\theta)\xi_\tau(r)^\top\big[b_{xx}(r,\bar{X}^{t,x;\bar{u}}(r)
                             +\theta\xi_\tau(r),\bar{u}(r))-\bar{b}_{xx}(r)\big]\xi_\tau(r)d\theta,\\
      \varepsilon_{\tau4}(r)&:=\int_0^1(1-\theta)\xi_\tau(r)^\top\big[\sigma_{xx}(r,\bar{X}^{t,x;\bar{u}}(r)+\theta\xi_\tau(r),\bar{u}(r))
                             -\bar{\sigma}_{xx}(r)\big]\xi_\tau(r)d\theta.
\end{aligned}
\right.
\end{equation*}
Similar to the proof of (\ref{Estimates of remainder terms of SDE}), using (\ref{estimate of SDEP in t}) we have, $\mathbf{P}\mbox{-}a.s.$,
\begin{equation}\label{Estimates of remainder terms of SDE-tau}
\left\{
\begin{aligned}
&\mathbb{E}\Big[\int_\tau^T|\varepsilon_{\tau1}(r)|^{2k}dr\big|\mathcal{F}_t^s\Big]\leq o(|\tau-s|^k),\ \mathbb{E}\Big[\int_\tau^T|\varepsilon_{\tau2}(r)|^{2k}dr\big|\mathcal{F}_t^s\Big]\leq o(|\tau-s|^k),\\
&\mathbb{E}\Big[\int_\tau^T|\varepsilon_{\tau3}(r)|^kdr\big|\mathcal{F}_t^s\Big](\omega)\leq o(|\tau-s|^k),\ \mathbb{E}\Big[\int_\tau^T|\varepsilon_{\tau4}(r)|^kdr\big|\mathcal{F}_t^s\Big](\omega)\leq o(|\tau-s|^k).
\end{aligned}
\right.
\end{equation}

Denote by $(Y^{\tau,\bar{X}^{t,x;\bar{u}}(s);\bar{u}}(\cdot),Z^{\tau,\bar{X}^{t,x;\bar{u}}(s);\bar{u}}(\cdot))$ the solution to the following BSDE (\ref{time disturbed controlled BSDE}) on $\big(\Omega,\mathcal{F},\{\mathcal{F}_\tau^t\}_{\tau\geq t},\mathbf{P}(\cdot|\mathcal{F}_\tau^t)\big)$ for $r\in[\tau,T]$:
\begin{equation}\label{time disturbed controlled BSDE}
\begin{aligned}
 &Y^{\tau,\bar{X}^{t,x;\bar{u}}(s);\bar{u}}(r)=\phi(X^{\tau,\bar{X}^{t,x;\bar{u}}(s);\bar{u}}(T))+\int_r^Tf\big(\alpha,X^{\tau,\bar{X}^{t,x;\bar{u}}(s);\bar{u}}(\alpha),\\
 &\quad Y^{\tau,\bar{X}^{t,x;\bar{u}}(s);\bar{u}}(\alpha),Z^{\tau,\bar{X}^{t,x;\bar{u}}(s);\bar{u}}(\alpha),\bar{u}(\alpha)\big)d\alpha
                                              -\int_r^TZ^{\tau,\bar{X}^{t,x;\bar{u}}(s);\bar{u}}(\alpha)dW(\alpha).
\end{aligned}
\end{equation}

For any $\tau\leq r\leq T$, set $\hat{Y}_\tau(r):=\big\langle p(r),\xi_\tau(r)\big\rangle+\frac{1}{2}\big\langle P(r)\xi_\tau(r),\xi_\tau(r)\big\rangle$. Applying It\^{o}'s formula, we get
\begin{equation}\label{applying Ito's formula-time}
d\hat{Y}_\tau(r)=C_\tau(r)dr+\hat{Z}_\tau(r)dW(r),\ r\in[\tau,T].
\end{equation}
where
\begin{equation*}
\left\{
{\small\begin{aligned}
             C_\tau(r)&:=-\big\langle\xi_\tau(r),\bar{f}_y(r)p(r)+\bar{f}_z(r)\bar{\sigma}_x^\top(r)p(r)
                       +\bar{f}_z(r)q(r)+\bar{f}_x(r)\big\rangle+\big\langle p(r),\varepsilon_{\tau3}(r)\big\rangle\\
                      &\quad+\big\langle q(r),\varepsilon_{\tau4}(r)\big\rangle+\mbox{tr}\Big\{-\frac{1}{2}\xi_\tau(r)\xi_\tau(r)^\top
                       \Big[\bar{f}_y(r)P(r)+2\bar{f}_z(r)\bar{\sigma}_x^\top(r)P(r)\\
                      &\quad+\bar{f}_z(r)Q(r)+\bar{\sigma}_{xx}^\top(r)\bar{f}_z(r)p(r)+\big[I_{n\times n},p(r),\bar{\sigma}_x^\top(r)p(r)+q(r)\big]\\
                      &\quad\cdot D^2\bar{f}(r)\big[I_{n\times n},p(r),\bar{\sigma}_x^\top(r)p(r)+q(r)\big]^\top\Big]
                       +\frac{1}{2}P(r)\varepsilon_{\tau5}(r)+\frac{1}{2}Q(r)\varepsilon_{\tau6}(r)\Big\},\\
       \hat{Z}_\tau(r)&:=\big\langle q(r),\xi_\tau(r)\big\rangle+\big\langle p(r),\bar{\sigma}_x^\top(r)\xi_\tau(r)
                       +\frac{1}{2}\xi_\tau(r)^\top\bar{\sigma}_{xx}(r)\xi_\tau(r)\big\rangle+\big\langle p(r),\varepsilon_{\tau4}(r)\big\rangle\\
                      &\quad+\frac{1}{2}\mbox{tr}\Big\{\xi_\tau(r)\xi_\tau(r)^\top Q(r)+P(r)\xi_\tau(r)\xi_\tau(r)^\top\bar{\sigma}_x^\top(r)+P(r)\bar{\sigma}_x(r)\xi_\tau(r)\xi_\tau(r)^\top\\
                      &\quad+P(r)\varepsilon_{\tau6}(r)\Big\},
\end{aligned}}
\right.
\end{equation*}
and
\begin{equation*}
\left\{
\begin{aligned}\varepsilon_{\tau5}(r)&:=\xi_\tau(r)\varepsilon_{\tau1}(r)^\top+\varepsilon_{\tau1}(r)\xi_\tau(r)^\top
                                      +\bar{\sigma}_x(r)\xi_\tau(r)\varepsilon_{\tau2}(r)^\top\\
                                     &\qquad+\varepsilon_{\tau2}(r)\xi_\tau(r)^\top\bar{\sigma}_x^\top(r)+\varepsilon_{\tau2}(r)\varepsilon_{\tau2}(r)^\top,\\
               \varepsilon_{\tau6}(r)&:=\xi_\tau(r)\varepsilon_{\tau2}(r)^\top+\varepsilon_{\tau2}(r)\xi_\tau(r)^\top.
\end{aligned}
\right.
\end{equation*}
Define $\widetilde{Y}_\tau(r):=Y^{\tau,\bar{X}^{t,x;\bar{u}}(s);\bar{u}}(r)-\hat{Y}_\tau(r),
\widetilde{Z}_\tau(r):=Z^{\tau,\bar{X}^{t,x;\bar{u}}(s);\bar{u}}(r)-\hat{Z}_\tau(r)$. Similar to the Steps 4, 5 in the proof of Theorem 3.1, we have
\begin{equation}\label{final estimate-time}
\begin{aligned}
\big|\widetilde{Y}_\tau(\tau)-\bar{Y}^{t,x;\bar{u}}(\tau)\big|^2&\leq o(|\tau-s|^2),\ \mathbf{P}\mbox{-}a.s.
\end{aligned}
\end{equation}

Note that $\big(\Omega,\mathcal{F},\mathbf{P}(\cdot|\mathcal{F}_\tau^t),W(\cdot)-W(\tau);u(\cdot))|_{[\tau,T]}\big)\in\mathcal{U}^w[\tau,T],\mathbf{P}\mbox{-}a.s.$ Thus by the definition of the value function $V$, we have
\begin{equation*}
V(\tau,\bar{X}^{t,x;\bar{u}}(s))\leq-Y^{\tau,\bar{X}^{t,x;\bar{u}}(s);\bar{u}}(\tau),\ \mathbf{P}\mbox{-}a.s.
\end{equation*}
Taking $\mathbb{E}(\cdot|\mathcal{F}_s^t)$ on both sides and noting that $\bar{X}^{t,x;\bar{u}}(s)$ is $\mathcal{F}_s^t$-measurable, we have
\begin{equation}\label{estimate inequality in t}
V(\tau,\bar{X}^{t,x;\bar{u}}(s))\leq\mathbb{E}\big[-Y^{\tau,\bar{X}^{t,x;\bar{u}}(s);\bar{u}}(\tau)\big|\mathcal{F}_s^t\big],\ \mathbf{P}\mbox{-}a.s.
\end{equation}
By the DPP of \cite{Peng97}, choose a common subset $\Omega_0\subseteq\Omega$ with $\mathbf{P}(\Omega_0)=1$ such that for any $\omega_0\in\Omega_0$, the following holds:
\begin{equation*}
\left\{
\begin{aligned}
&V(s,\bar{X}^{t,x;\bar{u}}(s,\omega_0))=-\bar{Y}^{t,x;\bar{u}}(s,\omega_0), \text { and }(\ref{estimate of SDEP in t}),(\ref{Estimates of remainder terms of SDE-tau}),(\ref{estimate inequality in t}) \mbox{ are satisfied for any }\\
&\mbox{rational } \tau>s,\big(\Omega,\mathcal{F},\mathbf{P}(\cdot|\mathcal{F}_s^t)(\omega_0),W(\cdot)-W(s);u(\cdot))|_{[s,T]}\big)\in\mathcal{U}^w[s,T], \mbox{ and }\\
&\sup\limits_{s\leq r\leq T}\big[|p(r,\omega_0)|+|P(r,\omega_0)|\big]<\infty.
\end{aligned}
\right.
\end{equation*}
Let $\omega_0\in\Omega_0$ be fixed, then for any rational number $\tau>s$, by (\ref{final estimate-time}) we have
\begin{equation}\label{final estimate-time-with remainder terms}
\begin{aligned}
    &V(\tau,\bar{X}^{t,x;\bar{u}}(s,\omega_0))-V(s,\bar{X}^{t,x;\bar{u}}(s,\omega_0))
    \leq\mathbb{E}\big[-Y^{\tau,\bar{X}^{t,x;\bar{u}}(s);\bar{u}}(\tau)+\bar{Y}^{t,x;\bar{u}}(s)\big|\mathcal{F}_s^t\big](\omega_0)\\
    &\leq\mathbb{E}\big[-Y^{\tau,\bar{X}^{t,x;\bar{u}}(s);\bar{u}}(\tau)+\bar{Y}^{t,x;\bar{u}}(\tau)-\bar{Y}^{t,x;\bar{u}}(\tau)
     +\bar{Y}^{t,x;\bar{u}}(s)\big|\mathcal{F}_s^t\big](\omega_0)\\
    &=\mathbb{E}\big[-\big\langle p(\tau),\xi_\tau(\tau)\big\rangle-\frac{1}{2}\big\langle P(\tau)\xi_\tau(\tau),\xi_\tau(\tau)\big\rangle\big|\mathcal{F}_s^t\big](\omega_0)+o(|\tau-s|)\\
    &\quad+\mathbb{E}\big[-\bar{Y}^{t,x;\bar{u}}(\tau)+\bar{Y}^{t,x;\bar{u}}(s)\big|\mathcal{F}_s^t\big](\omega_0)\\
    &=\mathbb{E}\Big[\int_s^\tau\bar{f}(r)dr-\big\langle p(\tau),\xi_\tau(\tau)\big\rangle-\frac{1}{2}\big\langle P(\tau)\xi_\tau(\tau),\xi_\tau(\tau)\big\rangle\Big|\mathcal{F}_s^t\Big](\omega_0)+o(|\tau-s|).
\end{aligned}
\end{equation}

Now let us estimate the terms on the right-hand side of (\ref{final estimate-time-with remainder terms}). To this end, we first note that for any $\varphi(\cdot),\hat{\varphi}(\cdot),\psi(\cdot)\in
L^2_\mathcal{F}([0,T];\mathbf{R}^n)$, we have
\begin{equation*}
\begin{aligned}
&\mathbb{E}\Big[\Big\langle\int_s^\tau\varphi(r)dr,\int_s^\tau\hat{\varphi}(r)dr\Big\rangle\Big|\mathcal{F}_s^t\Big](\omega_0)\\
&\leq\Big\{\mathbb{E}\Big[\big|\int_s^\tau\varphi(r)dr\big|^2\Big|\mathcal{F}_s^t\Big](\omega_0)\Big\}^{\frac{1}{2}}
 \Big\{\mathbb{E}\Big[\big|\int_s^\tau\hat{\varphi}(r)dr\big|^2\Big|\mathcal{F}_s^t\Big](\omega_0)\Big\}^{\frac{1}{2}}\\
&\leq(\tau-s)\Big\{\int_s^\tau\mathbb{E}\big[|\varphi(r)|^2\big|\mathcal{F}_s^t\big](\omega_0)dr
 \int_s^\tau\mathbb{E}\big[|\hat{\varphi}(r)|^2\big|\mathcal{F}_s^t\big](\omega_0)dr\Big\}^{\frac{1}{2}}=o(|\tau-s|),
\end{aligned}
\end{equation*}
as $\tau\downarrow s$, $\text{ for a.e. } s\in[t,T)$, and
\begin{equation*}
\begin{aligned}
&\mathbb{E}\Big[\Big\langle\int_s^\tau\varphi(r)dr,\int_s^\tau\psi(r)dW(r)\Big\rangle\Big|\mathcal{F}_s^t\Big](\omega_0)\\
&\leq\ \Big\{\mathbb{E}\Big[\big|\int_s^\tau\varphi(r)dr\big|^2\Big|\mathcal{F}_s^t\Big](\omega_0)\Big\}^{\frac{1}{2}}
 \Big\{\mathbb{E}\Big[\big|\int_s^\tau\psi(r)dW(r)\big|^2\Big|\mathcal{F}_s^t\Big](\omega_0)\Big\}^{\frac{1}{2}}\\
&\leq\ (\tau-s)^{\frac{1}{2}}\Big\{\int_s^\tau\mathbb{E}\big[|\varphi(r)|^2\big|\mathcal{F}_s^t\big](\omega_0)dr
 \int_s^\tau\mathbb{E}\big[|\psi(r)|^2\big|\mathcal{F}_s^t\big](\omega_0)dr\Big\}^{\frac{1}{2}}=o(|\tau-s|),
\end{aligned}
\end{equation*}
as $\tau\downarrow s$, $\text{ for a.e. } s\in[t,T)$. Each last equality in the above two inequalities holds, since the sets of right Lebesgue points possess full Lebesgue measures for integrable functions and $s\mapsto\mathcal{F}^t_s$ is right continuous in $s$. Thus by (\ref{first-order variational equation in t}) and (\ref{adjoint equation-first order}), we have
\begin{equation}\label{first part of difference of two value functions 2}
\begin{aligned}
 &\mathbb{E}\big[\langle p(\tau),\xi_\tau(\tau)\rangle\big|\mathcal{F}_s^t\big](\omega_0)
  =\mathbb{E}\Big[\langle p(s),\xi_\tau(\tau)\rangle+\langle p(\tau)-p(s),\xi_\tau(\tau)\rangle\big|\mathcal{F}_s^t\Big](\omega_0)\\
 &=\mathbb{E}\Big[\Big\langle p(s),-\int_s^\tau\bar{b}(r)dr-\int_s^\tau\bar{\sigma}(r)dW(r)\Big\rangle
  +\Big\langle-\int_s^\tau\Big(\bar{b}_x^\top(r)p(r)+\bar{f}_y(r)p(r)\\
 &\qquad+\bar{f}_z(r)\bar{\sigma}_x^\top(r)p(r)+\bar{f}_z(r)q(r)+\bar{\sigma}_x^\top(r)q(r)+\bar{f}_x(r)\Big)dr+\int_s^\tau q(r)dW(r),\\
 &\qquad-\int_s^\tau\bar{b}(r)dr-\int_s^\tau\bar{\sigma}(r)dW(r)\Big\rangle\Big|\mathcal{F}_s^t\Big](\omega_0)\\
 &=\mathbb{E}\Big[-\Big\langle p(s),\int_s^\tau\bar{b}(r)dr\Big\rangle
  -\int_s^\tau\big\langle q(r),\bar{\sigma}(r)\big\rangle dr\Big|\mathcal{F}_s^t\Big](\omega_0)+o(|\tau-s|).
\end{aligned}
\end{equation}
Similarly, by (\ref{second-order variational equation in t}) and (\ref{adjoint equation-second order}), we have
\begin{equation}\label{second part of difference of two value functions 2}
{\small\begin{aligned}
&\mathbb{E}\big[\xi_\tau(\tau)^\top P(\tau)\xi_\tau(\tau)\big|\mathcal{F}_s^t\big](\omega_0)
=\mathbb{E}\Big[\xi_\tau(\tau)^\top P(s)\xi_\tau(\tau)+\xi_\tau(\tau)^\top\big(P(\tau)-P(s)\big)\xi_\tau(\tau)\Big|\mathcal{F}_s^t\Big](\omega_0)\\
&=\mathbb{E}\Big[\Big(-\int_s^\tau\bar{b}(r)dr-\int_s^\tau\bar{\sigma}(r)dW(r)\Big)^\top P(s)
 \Big(-\int_s^\tau\bar{b}(r)dr-\int_s^\tau\bar{\sigma}(r)dW(r)\Big)\\
&\qquad+\Big(-\int_s^\tau\bar{b}(r)dr-\int_s^\tau\bar{\sigma}(r)dW(r)\Big)^\top\Big[-\int_s^\tau\Big(\big[\bar{f}_y(r)
 \bar{f}_z(s)\bar{\sigma}_x^\top(s)+\bar{b}_x^\top(s)\big]P(r)\\
&\qquad+P(r)\big[\bar{f}_z(r)\bar{\sigma}_x(r)+\bar{b}_x(r)\big]+\bar{\sigma}_x^\top(r)P(r)\bar{\sigma}_x(r)+\bar{f}_z(r)Q(r)+\bar{\sigma}_x^\top(r)Q(r)+Q(r)\bar{\sigma}_x(r)\\
&\qquad+\bar{b}_{xx}^\top(r)p(r)+\bar{\sigma}_{xx}^\top(r)\big[\bar{f}_z(r)p(r)+q(r)\big]+\big[I_{n\times n},p(r),\bar{\sigma}_x^\top(r)p(r)+q(r)\big]D^2\bar{f}(r)\\
&\qquad\cdot\big[I_{n\times n},p(r),\bar{\sigma}_x^\top(r)p(r)+q(r)\big]^\top\Big)dr+\int_s^\tau Q(s)dW(s)\Big]\Big(-\int_s^\tau\bar{b}(r)dr\\
&\qquad-\int_s^\tau\bar{\sigma}(r)dW(r)\Big)\Big|\mathcal{F}_s^t\Big](\omega_0)=\mathbb{E}\Big[\int_s^\tau\mbox{tr}\big\{\bar{\sigma}(r)^\top P(r)\bar{\sigma}(r)\big\}dr\Big|\mathcal{F}_s^t\Big](\omega_0)+o(|\tau-s|).
\end{aligned}}
\end{equation}
It follows from (\ref{final estimate-time-with remainder terms}), (\ref{first part of difference of two value functions 2}) and
(\ref{second part of difference of two value functions 2}), that for any rational $\tau\in(s,T]$ and at $\omega=\omega_0$,
\begin{equation}\label{difference of two value functions 4}
\begin{aligned}
    &V(\tau,\bar{X}^{t,x;\bar{u}}(s))-V(s,\bar{X}^{t,x;\bar{u}}(s))\leq\mathbb{E}\Big[\int_s^\tau\bar{f}(r)dr+\Big\langle p(s),\int_s^\tau\bar{b}(r)dr\Big\rangle\\
    & +\int_s^\tau\big\langle q(r),\bar{\sigma}(r)\big\rangle dr-\frac{1}{2}\int_s^\tau\mbox{tr}\big\{\bar{\sigma}(r)^\top P(r)\bar{\sigma}(r)\big\}dr\Big|\mathcal{F}_s^t\Big]+o(|\tau-s|)\\
   =&\ (\tau-s)\mathcal{H}_1(s,\bar{X}^{t,x;\bar{u}}(s),\bar{u}(s))+o(|\tau-s|).
\end{aligned}
\end{equation}
By definition
(\ref{right super- and sub-jets time}), we obtain that (\ref{connection-time}) holds for any (not only rational numbers) $\tau\in(s,T]$.
Finally, fix an $\omega\in\Omega$ such that (\ref{difference of two value functions 4}) holds for any $\tau\in(s,T]$. Then for any $\hat{q}\in D_{t+}^{1,-}V(s,\bar{X}^{t,x;\bar{u}}(s))$, by definition (\ref{right super- and sub-jets time}) and (\ref{difference of two value functions 4}) we have
\begin{equation*}
\begin{aligned}
0&\leq\liminf\limits_{\tau\downarrow s}\left\{\frac{V(\tau,\bar{X}^{t,x;\bar{u}}(s))-V(s,\bar{X}^{t,x;\bar{u}}(s))-\hat{q}(\tau-s)}{|\tau-s|}\right\}\\
 &\leq\liminf\limits_{\tau\downarrow s}\left\{\frac{(\mathcal{H}_1(s,\bar{X}^{t,x;\bar{u}}(s),\bar{u}(s))-\hat{q})(\tau-s)}{|\tau-s|}\right\}
\end{aligned}
\end{equation*}
Then, it is necessary that $\hat{q}\leq \mathcal{H}_1(s,\bar{X}^{t,x;\bar{u}}(s),\bar{u}(s))$.
Thus, (\ref{connection-time 2}) holds.
The proof is complete.\quad$\Box$

\section{Concluding Remarks}

This paper is the general extension of our companion paper \cite{NSW16}. In this paper, we have established a non-smooth version of the connection between the maximum principle and dynamic programming principle, for the stochastic recursive control problem when the control domain is non-convex and the diffusion coefficient depends on the control variable. By employing the viscosity solution of \cite{CIL92,YZ99}, the connection is now interpreted as four set inclusions. The first one is between the super-jet $D_x^{2,+}V(s,\bar{X}^{t,x;\bar{u}}(s))$ and $\{-p(s)\}\times[-P(s),\infty)$, the second one is between the sub-jet $D_x^{2,-}V(s,\bar{X}^{t,x;\bar{u}}(s))$ and $\{-p(s)\}\times(-\infty,-P(s)]$, the third one is between the right super-jet $D_{t+}^{1,+}V(s,\bar{X}^{t,x;\bar{u}}(s))$ and $[\mathcal{H}_1(s,\bar{X}^{t,x;\bar{u}}(s),\bar{u}(s)),\infty)$, and the fourth one is between the right sub-jet $D_{t+}^{1,-}V(s,\bar{X}^{t,x;\bar{u}}(s))$ and $(-\infty,\mathcal{H}_1(s,\bar{X}^{t,x;\bar{u}}(s),\bar{u}(s))]$. These new results have not only extended the classical one in \cite{Shi10} by eliminating the smoothness assumption on the value function, but also generalized the result obtained in \cite{NSW16} both to the second-order case and to the case with non-convex control domain.

Stochastic verification theorem of forward-backward controlled systems for viscosity solutions has been proved by \cite{Zhang12}. Based on this, we can derive the MP directly from DPP as the non-smooth version of Corollary 3.6. We have finished this topic and the result will appear elsewhere. Other extensions of the results in this paper to the jump-diffusion case (see Shi and Wu \cite{SW11}) and to the fully coupled forward-backward controlled stochastic systems (see Wu \cite{Wu97}, Yong \cite{Yong2010}), will be considered in the near future.

\vspace{2mm}

{\bf Acknowledgment.}\ The author thanks the associate editor and two anonymous referees for their careful reading of the previous version of the paper.

\end{document}